\documentclass[11pt,a4paper,twoside]{article}

\usepackage{amsmath}
\usepackage{amssymb}
\usepackage{amsthm}
\usepackage{graphicx}
\usepackage{psfrag}

\newcommand{\pa}[1]{{\partial#1}}

\newcommand{\h}[1]{{\hat#1}}
\newcommand{\ti}[1]{{\tilde#1}}

\newcommand{\ga}{\alpha}
\newcommand{\gb}{\beta}
\newcommand{\gd}{\delta}

\newcommand{\gga}{\gamma}
\newcommand{\gl}{\lambda}

\newcommand{\gr}{\rho}
\newcommand{\gs}{\sigma}
\newcommand{\gt}{\theta}
\newcommand{\gz}{\zeta}

\def\XXint#1#2#3{{\setbox0=\hbox{$#1{#2#3}{\int}$}
     \vcenter{\hbox{$#2#3$}}\kern-.5\wd0}}

\theoremstyle{plain}

\newtheorem{theorem}{Theorem}[section]

\newtheorem{remark}[theorem]{Remark}

\numberwithin{equation}{section} \numberwithin{figure}{section}

\begin{document}

\title{On the Linearization of the Painlev\'{e} III-VI Equations
and Reductions of the Three-Wave Resonant System}
\author{N.~Joshi,$^a$ A. V.~Kitaev,$^{ab}$ and P. A.~Treharne$^a$
\thanks{Electronic mail: nalini@maths.usyd.edu.au, kitaev@pdmi.ras.ru, philip@maths.usyd.edu.au}\\
$^a$School of Mathematics and Statistics, The University of Sydney,\\ NSW 2006, Australia,\\
$^b$Steklov Mathematical Institute, Fontanka 27, St Petersburg, 191023, Russia \\}
\date{June 12, 2007}
\maketitle

\abstract
We extend similarity reductions of the coupled (2+1)-dimensional three-wave resonant interaction
system to its Lax pair.
Thus we obtain new $3\times3$ matrix Fuchs--Garnier pairs for the third, fourth, and fifth Painlev\'e 
equations, together with the previously known Fuchs--Garnier pair for the sixth Painlev\'e equation.
These Fuchs--Garnier pairs have an important feature: they are linear with respect
to the spectral parameter. Therefore we can apply the Laplace transform to study these pairs.
In this way we found reductions of all pairs to the standard $2\times2$ matrix Fuchs--Garnier pairs obtained
by M. Jimbo and T. Miwa. As an application of the $3\times3$ matrix pairs, we found an integral auto-transformation
for the standard Fuchs--Garnier pair for the fifth Painlev\'e equation. It generates an Okamoto-like B\"acklund
transformation for the fifth Painlev\'e equation. Another application is an integral transformation relating two
different $2\times2$ matrix Fuchs--Garnier pairs for the third Painlev\'e equation.
\vspace{24pt}\vspace{12pt}\\
{\bf PACS 2006}: 02.30.Ik, 02.30.Gp, 02.30. Hq.\vspace{12pt}\vspace{24pt}\\
{\bf Running title}: On Linearization of the Painlev\'e III-VI Equations
\newpage
\setcounter{page}2

\section{Introduction}
 \label{sec:Intro}
The Painlev\'{e} equations are six classical nonlinear second-order
ordinary differential equations. They have been the subject of
intensive investigation in the last three decades, primarily due to
the fact that they appear in connection with a wide range of
physical problems, including soliton systems, quantum gravity,
string theory and random matrix theory. In this paper we will
concentrate on the third, fourth, fifth and sixth Painlev\'{e}
equations, the canonical forms of which are, respectively
\begin{align}
\text{P}_3:&\quad\frac{d^2y}{dt^2}=\frac1y\left(\frac{dy}{dt}\right)^2-\frac1t\frac{dy}{dt}+\frac1t
\Big(\ga y^2+\gb\Big)+\gga y^3+\frac{\gd}y,
 \label{P3}\\
\text{P}_4:&\quad\frac{d^2y}{dt^2}=\frac1{2y}\Big(\frac{dy}{dt}\Big)^2+\frac32y^3+4ty^2+2(t^2-\ga)y+\frac{\gb}y,
 \label{P4}\\
\text{P}_5:&\quad\frac{d^2y}{dt^2}=\Big(\frac1{2y}+\frac1{y-1}\Big)\Big(\frac{dy}{dt}\Big)^2-\frac1t\frac{dy}{dt}+
\frac{(y-1)^2}{t^2}\Big(\ga y+\frac{\gb}y\Big)+\frac{\gga y}t+\frac{\gd y(y+1)}{y-1},
 \label{P5}\\
\text{P}_6:&\quad\frac{d^2y}{dt^2}=\frac12\Big(\frac1y+\frac1{y-1}+\frac1{y-t}\Big)\Big(\frac{dy}{dt}\Big)^2-
\Big(\frac1t+\frac1{t-1}+\frac1{y-t}\Big)\frac{dy}{dt}+\notag\\
&\hspace{30mm}\frac{y(y-1)(y-t)}{t^2(t-1)^2}\Big(\ga+\frac{\gb t}{y^2}+\frac{\gga(t-1)}{(y-1)^2}+
\frac{\gd t(t-1)}{(y-t)^2}\Big),
 \label{P6}
\end{align}
where $\ga$, $\gb$, $\gga$, and $\gd$ are arbitrary complex parameters,
see \cite{I1956}.

In 1888 L. Fuchs found that if the second order ODE,
\begin{equation}
 \label{eq:INT-L.Fuchs1}
\frac{d^2\psi}{dx^2}=p(x)\psi,
\end{equation}
where $p(x)$ is a rational function, has a monodromy group independent of
the position of singular points, $x = t_i$, $i=1,\ldots,n$
({\it isomonodromy deformation}), then the function $\psi$ satisfies one more auxiliary equation,
\begin{equation}
 \label{eq:INT-L.Fuchs2}
\frac{\partial\psi}{\partial t_i}= A_i(x)\psi +B_i(x)\frac{\partial\psi}{\partial x},
\end{equation}
where $A_i$ and $B_i$ are some rational functions of $x$. In 1905, R. Fuchs reported that in the simplest
nontrivial particular case when the equation is of the Fuchsian type with four singular points located at
$0,1,t_1=t,\infty$ and the fifth singular point is an apparent singularity, then its location, $y(t)$, is
governed by $P_6$  (see extended later article by R. Fuchs~\cite{F1907}). A few years later R. Garnier~\cite{G1912}
considered the general case of the Fuchsian equation \eqref{eq:INT-L.Fuchs1} of the second order and derived the
generalization of $P_6$ which is now known as the Garnier system. In the same paper he also found the
pairs \eqref{eq:INT-L.Fuchs1}, \eqref{eq:INT-L.Fuchs2} for the other Painlev\'e equations. In the latter case,
equation \eqref{eq:INT-L.Fuchs1} is non-Fuchsian. In honor of this contribution, we call the pairs that define
isomonodromy deformations of linear ODEs (of arbitrary order) with rational coefficients, Fuchs--Garnier pairs.

The Fuchs--Garnier pairs associated with each Painlev\'e equation play a very important role in the theory and
applications of the Painlev\'e equations. Nowadays, as a result of the intensive studies of the Painlev\'e equations,
many different Fuchs--Garnier pairs have been derived
\cite{B2005a,B2005b,FN1980,H1994,JM1981II,KH1999,K1985,KV,K2006,MCB1997,NY1999,NY2000,WH2003}.
The methods and ideas used in these derivations vary widely. As a result, most of the Painlev\'e equations possess
a few different Fuchs--Garnier pairs whose equivalence is not yet established. Such Fuchs--Garnier pairs very often
contain matrix differential equations. Thus, 
different Fuchs--Garnier pairs for the same Painlev\'e equation can have
different matrix dimension, different analytic structure and even, if the two first features are the same, they can
still have different {\it parametrization} of the matrix elements by the Painlev\'e functions. Our general belief is
that all different Fuchs--Garnier pairs associated with the same Painlev\'e equation
should be related by some explicit transformations.

These transformations are interesting not only from the purely theoretical point of
view, but also from a practical one. For example, for the scalar and $2\times2$ matrix equations the corresponding
analytic and asymptotic theories are much simpler and better developed than those for the multidimensional cases,
therefore it might be useful to transport results obtained for the scalar and $2\times2$ matrix Fuchs--Garnier 
pairs to the multidimensional case. Also, in applications, such as in geometry, the solutions of the Fuchs--Garnier 
pairs often have a very definite interpretation, e.g., as the functions defining embeddings of some surfaces. 
Therefore, explicit relations between different Fuchs--Garnier pairs, even with the same matrix dimensions, might 
lead to interesting insights in geometry and mathematical physics.

It is clear from the definition given above that the role of the two equations in each Fuchs--Garnier pair is not
symmetric. There is a \lq\lq defining" equation, namely equation~\eqref{eq:INT-L.Fuchs1} and the \lq\lq deformation"
equation, namely equation~\eqref{eq:INT-L.Fuchs2}. The independent variable of the defining equation is called
a {\it spectral} variable (parameter); we denote it by $x$ or $\lambda$. The coefficients of both equations in
Fuchs--Garnier pairs are rational functions of this variable. When the defining equation is given, the deformation
equation can be derived from the isomonodromy condition. Therefore, sometimes for brevity, to present the
Fuchs--Garnier pair, we write only one defining equation.

Together with the original scalar Fuchs--Garnier pairs, the $2\times2$ matrix versions first presented by M. Jimbo
and T. Miwa \cite{JM1981II} play an important role in the study of the Painlev\'e equations. The defining equation
\begin{equation} \label{JMlam}
\frac{dY}{dx} = A^{n}(x;t) Y,
\end{equation}
has the following particular forms for the Painlev\'e equations $P_n$ listed above:
\begin{subequations}
 \label{JMlam2}
\begin{align}
A^{3}(x;t) &= \frac{A^{3}_{0}(t)}{x^2} + \frac{A^{3}_{1}(t)}{x} + A^{3}_{2}(t), \label{P3-JMlam} \\
A^{4}(x;t) &= \frac{A^{4}_{0}(t)}{x} + A^{4}_{1}(t) + x A^{4}_{2}(t), \label{P4-JMlam} \\
A^{5}(x;t) &= \frac{A^{5}_{0}(t) }{x} + \frac{A^{5}_{1}(t)}{x-1} + A^{5}_{2}(t), \label{P5-JMlam} \\
A^{6}(x;t) &= \frac{A^{6}_{0}(t)}{x} + \frac{A^{6}_{1}(t)}{x-1}+\frac{A^{6}_{t}(t)}{x-t}
\label{P6-JMlam}
\end{align}
\end{subequations}
The matrices $A^k_i(t)$ are independent of the spectral parameter and are parameterized by the solutions
of the corresponding Painlev\'e equations (see Appendix C in \cite{JM1981II}). To distinguish other
$2\times2$ matrix Fuchs--Garnier pairs that are known for the same Painlev\'e equations, we call these
pairs Fuchs--Garnier pairs in the {\it Jimbo--Miwa parametrization}. For convenience of the reader,
we present the Fuchs--Garnier pair for $P_5$ in Jimbo--Miwa parametrization in Appendix~\ref{app:A}.

In his studies of the Painlev\'e equations, Okamoto pointed out that the Painlev\'e equations have subgroups
of symmetries isomorphic to some affine Weyl groups, \cite{O1987a}--\cite{O1986}. Using this fact, he constructed
nonlinear representations of these groups as birational canonical transformations of the Hamiltonian systems
associated with the Painlev\'e equations. As we explain in Appendix~\ref{app:A} on the example of $P_5$, there
is a problem with finding the linear representation of these affine Weyl groups in the space of solutions of the
Fuchs--Garnier pairs in the Jimbo--Miwa parametrization. Since the latter pairs proved to be a highly effective
and convenient tool for the complete description of global asymptotic properties of all solutions of the Painlev\'{e}
equations and in various applications, there is a motivation to complete the theory of these Fuchs--Garnier pairs
with the representation of the affine Weyl symmetries.

The goal of this paper is to create useful tools for finding transformations of solutions of
Fuchs --Garnier pairs that can answer the questions raised above. Our main stimulus in this work was a
recent understanding of some of the questions raised above for the case of $P_6$ in the work by 
M. Mazzocco~\cite{M2002}, D. Novikov~\cite{N2006} and G. Filipuk~\cite{F2006}. The latter two works explain that the linear representation
of one nontrivial, from the isomonodromy point of view, case of Okamoto's affine Weyl symmetries for $P_6$ is given
by the Euler integral auto-transform for the Fuchs--Garnier pair in the Jimbo--Miwa parametrization.
The work by Mazzocco explains that the \lq\lq dual" Fuchs--Garnier pair for $P_6$ found by J. Harnad~\cite{H1994}
can be mapped to the Fuchs--Garnier pair in Jimbo--Miwa parametrization by the Laplace transform\footnote{
The fact that defining the equation of Harnad's Fuchs--Garnier pair is related with the Jimbo--Miwa one via the
Laplace transform was observed by Dubrovin \cite{D1999}. In this connection, the work by
W. Balser, W. Jurkat, and D. Lutz \cite{BJL1981} should be mentioned also.}.
Moreover, the linear representation of the Okamoto transformation for $P_6$ is just a multiplication of the solution
by the scalar factor $\lambda^\alpha$ for a suitable choice of the parameter $\alpha$. Since multiplication by
$\lambda^\alpha$ is conjugate by the Laplace transform to the Euler transformation the result of the
works \cite{N2006} and \cite{F2006} follows immediately. So, the $3\times3$ Fuchs--Garnier pair by Harnad
serves as a useful auxiliary object in this study with the Laplace transform as the main instrument.

Our objective is to extend the ideas related with the Laplace transform to the other Painlev\'e equations. For
this purpose, we have to find proper analogues of Harnad's Fuchs-Garnier pair for the other Painlev\'e equations.
The adjective \lq\lq proper" here means that the pairs should be $3\times3$ matrix equations and we should be able
to apply to them the Laplace transform. The latter condition suggests that at least the defining member of the
Fuchs--Garnier pair, the ODE with respect to the spectral parameter, $\lambda$, should have linear coefficients in
$\lambda$, i.e.,
\begin{equation}
 \label{eq:INT-secondary-linearization}
\big(\lambda B_1(t)+B_2(t)\big)\frac{d\Phi}{d\lambda}=\Big(\lambda B_3(t)+B_4(t)\Big)\Phi.
\end{equation}
We note that in the work by J. Harnad mentioned above, the Fuchs--Garnier pair with \lq\lq spectral equation"
\eqref{eq:INT-secondary-linearization} where $B_2(t)=0$, $\det B_1\neq0$ was found for $P_6$.
Because of that result, it is easy to understand that Fuchs--Garnier pairs with the \lq\lq defining equation"
\eqref{eq:INT-secondary-linearization} in $3\times3$ matrices should exist for all Painlev\'e
equations. Actually, M. Noumi and Y. Yamada~\cite{NY1999, NY2000} found such a pair for the symmetric version 
of $P_4$. The latter pair was further studied by A. Sen, A. Hone, and P.A. Clarkson, however, in these studies,
the Laplace transform was not applied and the relation with the Jimbo--Miwa Fuchs--Garnier
pair \eqref {JMlam}, \eqref{P4-JMlam} was not yet realized. In this paper, we report
the Fuchs--Garnier pairs of the type~\eqref{eq:INT-secondary-linearization} for $P_3-P_6$: the pairs for
$P_3$, $P_4$, and $P_5$ are new, the pair for $P_6$ coincides with the known one \cite{H1994,M2002}.
We note that our Fuchs--Garnier pair for $P_4$ has a different singularity structure comparing to the one by
Noumi--Yamada~\cite{NY2000}. However, as we show, there is an invertible integral transformation linking together
these pairs. We also establish the relation of the new pairs for $P_3-P_5$ and the Noumi--Yamada pair to
the Jimbo--Miwa Fuchs--Garnier pairs \eqref{JMlam}, \eqref{P3-JMlam}--\eqref{P5-JMlam}, together with the known
result for $P_6$.

Let us remark that the phrase \lq\lq linearization of the Painlev\'e equations" is widely understood to mean an
association with some Fuchs--Garnier pair. In this paper, we extend this phrase to a \lq\lq secondary" linearization,
i.e., association of the Painlev\'e equations with Fuchs--Garnier pairs of the
form~\eqref{eq:INT-secondary-linearization}. Looking ahead, we note that secondary
linearization\footnote{In matrices with the dimension higher than 3, of course.} is possible for any so-called
higher-order Painlev\'e equations, however that is already the subject of another story.

The general form \eqref{eq:INT-secondary-linearization} for the defining equation of the 
Fuchs--Garnier pairs has matrix dimension three and, in the general case, contains too many variables for 
linear representations of the Painlev\'{e} equations. Instead of analyzing the general case as our starting 
point, we chose another and faster way to find the the proper Fuchs--Garnier pairs, which is based on the 
following observations.

Recently, there appeared two independent works by R.~Conte, A.~M.~Grundland, and
M.~Musette \cite{CGM2006} and S.~Kakei, T.~Kikuchi \cite{KK2007} were the authors obtain Harnad's Fuchs-Garnier
pair for $P_6$ by using an extension of the similarity reduction \cite{K1990} for the three-wave resonant
interaction (3WRI) system in $(1+1)$ (one spatial and one time) dimensions, to the corresponding Lax pair. The
Lax pair for this system was given in terms of two commuting first order differential operators in $3\times3$
matrices by V. E. Zakharov and S. V. Manakov \cite{ZM1973}.
The authors of  \cite{CGM2006} and \cite{KK2007} were able to get the Fuchs--Garnier pair already studied by Harnad
and Mazzocco and used their parametrization to get explicit formulae for the solutions of 3WRI system in terms of
$P_6$ with the complete set of the coefficients.

L. Martina and P. Winternitz in \cite{MW1989} obtained all classical similarity reductions for $(2+1)$ 3WRI system:
\begin{equation} \label{3WRIcharac}
\frac{\pa{u_{j}}}{\pa{x_{j}}} = i u_{m}^{*} u_{n}^{*}, \quad
\frac{\pa{u^{*}_{j}}}{\pa{x_{j}}} = -i u_{m} u_{n}, \quad i^{2} =
-1,
\end{equation}
where $(j,m,n)$ denotes any cyclic permutation of $(1,2,3)$, $u_{j},
u_{j}^{*}$ are the complex amplitudes of the wave packets, and star
denotes complex conjugation. This system is also integrable and it possesses, of course, a much richer group of
symmetries and corresponding similarity reductions than that in $(1+1)$ dimensions.
In particular, Martina and Winternitz found reductions
to the $P_6$, $P_5$, $P_4$, and $P_3$ equations with the complete set of the coefficients. It is important to note
that {\it the reductions to $P_3$, $P_4$, and $P_5$  cannot be restricted to the $(1+1)$ case of \emph{3WRI} system}.

For each reduction Martina and Winternitz used group theoretical methods to reduce the system of three complex PDEs
\eqref{3WRIcharac} to a system of three complex ODEs of the first order. In the latter system, they separated real
and imaginary part to arrive at a system of six real ODEs of the first order. They showed that for all similarity
reductions three of the six ODEs can be converted to one ODE of the third order which possesses the Painlev\'e
property while the rest three can be solved in quadratures in terms of the solution of the third order equation.
Such third order equations can always be integrated once to give a quite complicated ODE quadratic with respect to
the second derivative, a so-called SD equation (second order second degree  ODE). The latter ODEs was integrated
by Bureau \emph{et al} \cite{B1972,BGG1972} in terms of the solutions of the Painlev\'e equations mentioned
above\footnote{For recent advanced results see the paper by C.~M.~Cosgrove and G.~Scoufis \cite{CS1993}.}.
The similarity solutions obtained in this way in most cases are not explicitly written in terms of the canonical
Painlev\'e functions, because in the papers \cite{B1972,BGG1972} solutions of the SD equations are not always
presented in simple form in terms of the canonical functions. So in this study no any techniques related with the
Lax pairs were involved.

We note that solutions of the 3WRI system~\eqref{3WRIcharac} are not analytic, and therefore working with that
system we cannot achieve our goal -- to get Fuchs--Garnier pairs for the general case of the Painlev\'e
equations -- without any artificial restrictions. So, we have to consider an analytic extension of the 3WRI system
which we call also the coupled 3WRI system. We do it in the standard way, namely, we forget that the upper script 
$^*$ means complex conjugation in the six equations in \eqref{3WRIcharac} and we consider $u_j$ and $u_j^*$ as 
independent complex functions. The coupling procedure spoils neither its integrability, so that formally the same 
Lax pair serves for the coupled version of 3WRI system, nor the Martina--Winternitz similarity reductions.

To construct secondary linearized Fuchs--Garnier pairs for the Painlev\'e equations~\eqref{P3}--\eqref{P6},
we have to find for each Martina--Winternitz similarity reduction its extension to the Lax pair for the coupled
3WRI system. At this stage we arrive at $3\times3$ matrix Fuchs--Garnier pairs for a system of ODEs defining
similarity solutions of the coupled 3WRI system. A substantial question here is how to introduce the spectral
parameter; since originally the Lax pair for the 3WRI system does not contain any spectral parameters,
this is the major difference with the situation for $(1+1)$ integrable systems, where the Lax pairs already possess
the spectral parameter. While there are many papers in the literature concerning extensions of the similarity
reductions of $(1+1)$ integrable systems to their Lax pairs, this methodology is well known since the work of
H. Flashka and A. C. Newell \cite{FN1980}, we do not know such works for $(2+1)$ integrable systems. We note that
in our case the reduction cannot be done successively: $(2+1)\to(1+1)\to(1+0)$.
We found that for all similarity reductions it is possible to introduce the spectral parameter such that the 
defining equations of the resulting Fuchs--Garnier pairs gain the form \eqref{eq:INT-secondary-linearization}.

After the $3\times3$ matrix Fuchs--Garnier pairs are obtained we use Laplace and/or gauge transformations to
map these pairs to the $2\times2$ Fuchs--Garnier pairs in the Jimbo--Miwa parametrization. Comparing
parameterizations between the one that comes from the coupled 3WRI system and the Jimbo--Miwa parametrization
we obtain explicit formulae for the similarity solutions in terms of the canonical Painlev\'e functions. This
comparison also allows us to parameterize the $3\times3$ Fuchs--Garnier pairs in terms of the canonical
Painlev\'e functions, i.e., to obtain secondary linearized Fuchs--Garnier pairs for the Painlev\'e equations.

We also consider parametrization of the similarity solutions for the physical case of our coupled 3WRI system,
i.e., the original 3WRI system. At this stage we also arrive at SD functions but in this approach they have a
lucid sense as the Hamiltonians ($\tau$-functions) for the Painlev\'e equations.

The paper consists of six Sections and one Appendix. Section~\ref{sec:Intro} is the Introduction.
In Section~\ref{sec:3WRI} we recall the Lax pair for the 3WRI system. Sections~\ref{sec:P6}--\ref{sec:P3}
represent the main body of the paper: each one is devoted to the corresponding Painlev\'e equation beginning with
$P_6$ and finishing with $P_3$. The Sections are divided into Subsections which represent logical steps of the
derivation indicated above: similarity reductions from Martina and Winternitz, extensions of the reduction to the
Lax pair, reductions via the Laplace transform to the Fuchs--Garnier pairs in the Jimbo--Miwa parameterizations,
parameterizations of similarity solutions by the Painlev\'e functions. Sections~\ref{sec:P5} and \ref{sec:P3} have
also Subsections with the alternate reductions of the $3\times3$ matrix Fuchs--Garnier pairs to the $2\times2$ ones.
Section \ref{sec:P5} contains also one more extra subsection with the derivation of the Okamoto transformation.
Appendix A is devoted to the spectral interpretation of the B\"acklund transformations for $P_5$. In particular,
we define the Okamoto transformation and, at the very end, present the alternate parametrization of isomonodromy
deformations for equation~\eqref{JMlam}, \eqref{P5-JMlam}.

The main results obtained in this work are as follows:
\begin{enumerate}
\item
We introduced a notion of secondary linearization for the Painlev\'e equations as
the Fuchs--Garnier pairs with the defining equation~\eqref{eq:INT-secondary-linearization} in $3\times3$
matrices. We found these pairs for $P_3$--$P_6$. Three pairs for $P_3$, $P_4$, and $P_5$ are new. The pair for
$P_6$ coincides with the Harnad one~\cite{H1994}, see Subsections~\ref{subsec:P6-FG-Laplace},
\ref{subsec:P5-FG-Laplace}, \ref{subsec:P4-FG-Laplace}, \ref{subsec:P3-FG-Laplace}.
We prove that the pair for $P_4$ found in this paper is equivalent to the pair for the symmetric form of
this equation by Noumi and Yamada~\cite{NY1999,NY2000}, but the singularity structure of our one is different,
see Subsection~\ref{subsec:P4alt};
\item
We found a relation of all secondary linearized pairs to the Fuchs--Garnier pairs in the Jimbo--Miwa
parametrization. This is a new result only for the pairs for $P_3$, $P_4$, and $P_5$:
see Subsections~\ref{subsec:P5-FG-Laplace}, \ref{subsec:P4-FG-Laplace}, \ref{subsec:P3-FG-Laplace};
\item
For $P_5$ we found an explicit linear representation for the nontrivial Okamoto affine Weyl symmetry.
It is given as an integral auto-transform of the solution $Y$ of the Jimbo--Miwa Fuchs--Garnier pair.
The mechanism of its derivation is different from that for the analogous result for $P_6$,
see Subsection~\ref{subsec:Okamoto-BT};
\item
For the $2\times2$ Fuchs--Garnier pair for $P_5$ with the defining equation~\eqref{JMlam}, \eqref{P5-JMlam}
we found a simpler and more natural parametrization, which we call the \lq\lq true" Jimbo--Miwa parametrization,
see Appendix~\ref{app:A};
\item
For both cases of $P_3$ (the complete and degenerate) two $2\times2$ Fuchs--Garnier pairs are known
see, e.g. \cite{KH1999,K2006}, we found that they are related via an integral transform, see
Subsection~\ref{subsec:P3alt};
\item
As a byproduct of our work,
for both coupled and physical cases of 3WRI system and for all similarity reductions to the Painlev\'e
equations explicit parameterizations in terms of the canonical Painlev\'e functions are obtained:
see Subsections~\ref{subsec:P6-parametrization}, \ref{subsec:P5-parametrization}, \ref{subsec:P4-parametrization},
\ref{subsec:P3-parametrization}.
\end{enumerate}

The secondary linearization also exists, of course, for the first and second Painlev\'e equations.
They are not related with the similarity reductions of 3WRI system and corresponding results will be published
separately.

We expect that this approach with the auxiliary secondary linearized Fuchs--Garnier pairs will be
very fruitful for the hierarchies of the Painlev\'e equations.


\section{Lax Pair for the 3WRI System}
 \label{sec:3WRI}

System \eqref{3WRIcharac} admits a Lax pair found by Kaup
\cite{K1980}. We write it here in a modified form with the spectral
parameter $k$:
\begin{equation} \label{3WRI-Lax}
\begin{split}
\frac{\pa{\psi}_{j}}{\pa{x}_{m}} - ik \kappa_{m} \psi_{j} &= -i u_{n}^{*} \psi_{m} \\
\frac{\pa{\psi}_{m}}{\pa{x}_{j}} - ik \kappa_{j} \psi_{m} &= i u_{n}
\psi_{j}
\end{split}
\end{equation}
where $(j,m,n)$ denotes any cyclic permutation of $(1,2,3)$,
$\psi_{j} = \psi_{j}(x_{m},k)$ are scalar functions, $\kappa_{j}$
are real constants, and $k \in \mathbb{C}$ is the spectral
parameter. We note that our notation differs from Kaup's one by the
factor, $\psi_{j} \mapsto
\psi_{j}\exp{[ik(\kappa_{1}x_{1}+\kappa_{2}x_{2}+\kappa_{3}x_{3})]}$,
for all $j=1,2,3$. The spectral parameter $k$ appeared in Kaup's
analysis of the scattering problem for the system \eqref{3WRI-Lax}
for a different class of solutions of \eqref{3WRIcharac}.

System \eqref{3WRI-Lax} can be written in matrix form in the
following way
\begin{equation}\label{3WRI-LP}
\begin{split}
\mathcal{D}_{1} \Psi &= i\big( k K_{1} + U_{1} \big) \Psi \\
\mathcal{D}_{2} \Psi &= i\big( k K_{2} + U_{2} \big) \Psi
\end{split}
\end{equation}
where $\Psi$ is a $3\times3$ matrix-valued function, the matrix
operators $\mathcal{D}_{1}, \mathcal{D}_{2}$ and the matrices
$K_{1}, K_{2}$ and $U_{1}, U_{2}$ are defined as follows:
$$
\mathcal{D}_{1} = \mathrm{diag}[ \pa{}_{x_{2}}, \pa{}_{x_{3}},
\pa{}_{x_{1}} ],\qquad \mathcal{D}_{2} = \mathrm{diag}[
\pa{}_{x_{3}}, \pa{}_{x_{1}}, \pa{}_{x_{2}} ],
$$
\begin{align*}
K_{1} &=
\begin{pmatrix}
\kappa_{2} & 0 & 0 \\
0 & \kappa_{3} & 0 \\
0 & 0 & \kappa_{1}
\end{pmatrix}, \qquad
U_{1} =
\begin{pmatrix}
0 & -u_{3}^{*} & 0 \\
0 & 0 & -u_{1}^{*} \\
-u_{2}^{*} & 0 & 0
\end{pmatrix}, \\
K_{2} &=
\begin{pmatrix}
\kappa_{3} & 0 & 0 \\
0 & \kappa_{1} & 0 \\
0 & 0 & \kappa_{2}
\end{pmatrix}, \qquad
U_{2} =
\begin{pmatrix}
0 & 0 & u_{2} \\
u_{3} & 0 & 0 \\
0 & u_{1} & 0
\end{pmatrix}.
\end{align*}
We note that, when written in standard cartesian coordinates, the
linear system \eqref{3WRI-LP} is equivalent to the Lax pair
considered by Fokas and Ablowitz in \cite{FA1984}.

In the following sections we will investigate the particular
similarity reductions found in \cite{MW1989} that are linked to $P_3$
$P_4$, $P_5$ and $P_6$, giving an explicit extension of each
reduction on the Lax pair \eqref{3WRI-LP}.


\section{Similarity Reduction to the Sixth Painlev\'{e} Equation}
 \label{sec:P6}

The following similarity reduction for system \eqref{3WRIcharac} was
obtained in \cite{K1990} and \cite{MW1989}:
\begin{equation} \label{SimRedVI}
u_{1} = (x_{2}-x_{3})^{-1+i\gr_{1}} v_{1}, \quad u_{2} =
(x_{1}-x_{3})^{-1+i\gr_{2}} v_{2}, \quad u_{3} =
(x_{1}-x_{2})^{-1+i\gr_{3}} v_{3},
\end{equation}
where $v_j=v_{j}(\tau)$ with
\begin{equation}
\tau = \frac{x_{1}-x_{3}}{x_{2}-x_{3}},
\end{equation}
and $\gr_{1}, \gr_{2}, \gr_{3}$ are real constants such that
\begin{equation}
 \label{eq:P6-rho-relation}
\gr_1+\gr_2+\gr_3=0.
\end{equation}
Under this reduction, system~\eqref{3WRIcharac} reduces to the following system of ODEs:
\begin{equation} \label{RedSysVIa}
\begin{aligned}
\tau^{1+i\gr_{2}} (\tau-1)^{1+i\gr_{3}} v_{1}' &= iv_{2}^{*}v_{3}^{*}\\ 
\tau^{i\gr_{2}} (\tau-1)^{1+i\gr_{3}} v_{2}' &= -iv_{3}^{*}v_{1}^{*}\\ 
\tau^{1+i\gr_{2}} (\tau-1)^{i\gr_{3}} v_{3}' &= iv_{1}^{*}v_{2}^{*},
\end{aligned}
\end{equation}
where prime denotes differentiation with respect to $\tau$.

It is mentioned in the Introduction that the above system was integrated directly
in \cite{MW1989} in terms of the general solution of an SD equation which,
in turn, is solvable in terms of the sixth Painlev\'{e} function.
We will show at the end of Subsection~\ref{subsec:P6-parametrization}
that the similarity solutions can be written in a (relatively) simple way in terms of
the canonical $P_6$ functions, so that in this case the $SD$ equation is an
intermediate object that makes the formulae cumbersome. For the other similarity
reductions, SD functions are actually needed.

The one-dimensional restriction of the similarity reduction
\eqref{SimRedVI} was used in the recent works \cite{CGM2006} and
\cite{KK2007} to obtain a $3\times3$ Fuchs--Garnier pair for $P_6$
from the (1+1)-dimensional scattering Lax pair. In the remainder of this
section we will rederive this result from the (2+1)-dimensional
perspective.
\begin{remark}{\rm
However before we generalize this similarity reduction to the
coupled case of the 3WRI system. One adds to \eqref{SimRedVI} and
\eqref{RedSysVIa} the formally conjugated equations
\begin{equation*}
u_1^*=(x_2-x_3)^{-1-i\gr_1}v_1^*,\quad
u_2^*=(x_1-x_3)^{-1-i\gr_2}v_2^*,\quad
u_3^*=(x_1-x_2)^{-1-i\gr_3}v_3^*,\\
\end{equation*}
\begin{equation}
 \label{RedSysVIb}
\begin{aligned}
\tau^{1-i\gr_2}(\tau-1)^{1-i\gr_3}{v_1^*}'&=-iv_2v_3,\\
\tau^{-i\gr_2}(\tau-1)^{1-i\gr_3}{v_2^*}'&=iv_3v_1,\\
\tau^{1-i\gr_2}(\tau-1)^{-i\gr_3}{v_3^*}'&=-iv_1v_2.
\end{aligned}
\end{equation}
Note that in the coupled case $\rho_1,\rho_2,\rho_3\in\mathbb C$ satisfy
the same relation \eqref{eq:P6-rho-relation} and the functions
$v_j$ and $v_j^*$ are not assumed to be complex conjugates. In the
most part of this Section we deal with the coupled 3WRI system and
turn back to the physical case at the end of Subsection~\ref{subsec:P6-parametrization}.}
\end{remark}

\subsection{The $3\times3$ Fuchs--Garnier Pair}

To compute the reduced Lax pair we introduce the spectral parameter
$\gl$ in the following way
\begin{equation*} \label{lambdaVI}
\gl = (x_{2}-x_{3}) k.
\end{equation*}
Writing $\Psi(x_{j},k) = R(x_{j})\ti{\Phi}(\tau,\gl)$, where
$R(x_{j})$ is given by
\begin{subequations}
\begin{equation*} \label{R-6}
R(x_{1},x_{2},x_{3}) = \mathrm{diag}\, \big(
(x_{2}-x_{3})^{i\gt_{23}}, (x_{2}-x_{3})^{i\gt_{31}},
(x_{2}-x_{3})^{i\gt_{12}} \big),
\end{equation*}
and
\begin{equation*} \label{theta-6}
\gt_{12} - \gt_{31} = \gr_{1}, \quad \gt_{23}-\gt_{12} = \gr_{2},
\quad \gt_{31}-\gt_{23} = \gr_{3},
\end{equation*}
\end{subequations}
we find that Lax pair \eqref{3WRI-LP} can be rewritten as follows:
\begin{align*}
C_{1}\ti{\Phi}_{\tau} + \gl D_{1} \ti{\Phi}_{\gl} &= i\big( \gl K_{1} + V_{1} \big) \ti{\Phi} \\
C_{2}\ti{\Phi}_{\tau} + \gl D_{2} \ti{\Phi}_{\gl} &= i\big( \gl
K_{2} + V_{2} \big) \ti{\Phi},
\end{align*}
where the matrices $C_{j}, D_{j}, K_{j}, V_{j}$ are given by
\begin{align*}
C_{1} &= \mathrm{diag}\, \big( -\tau, \tau-1, 1 \big),& C_{2} &= \mathrm{diag}\, \big( \tau-1, 1, -\tau \big), \\
D_{1} &= \mathrm{diag}\, \big( 1, -1, 0 \big),& D_{2} &= \mathrm{diag}\, \big( -1, 0, 1 \big), \\
K_{1} &= \mathrm{diag}\, \big( \kappa_{2},\kappa_{3},\kappa_{1} \big),& K_{2}
&= \mathrm{diag}\, \big( \kappa_{3},\kappa_{1},\kappa_{2} \big), \\
V_{1} &=
\begin{pmatrix}
-\gt_{23} & -(\tau-1)^{-1-i\gr_{3}}v_{3}^{*} & 0 \\
0 & \gt_{31} & -v_{1}^{*} \\
-\tau^{-1-i\gr_{2}}v_{2}^{*} & 0 & 0
\end{pmatrix},& & \\
V_{2} &=
\begin{pmatrix}
\gt_{23} & 0 & \tau^{-1+i\gr_{2}}v_{2} \\
(\tau-1)^{-1+i\gr_{3}}v_{3} & 0 & 0 \\
0 & v_{1} & -\gt_{12}
\end{pmatrix}. & &
\end{align*}
After rearranging the above system, we find
\begin{subequations} \label{LP-VIgen}
\begin{align}
\ti{\Phi}_{\gl} &= \Big( Q^{(1)} + \frac{Q^{(0)}}{\gl} \Big) \ti{\Phi} \label{LP-VI1} \\
\ti{\Phi}_{\tau} &= \Big( \gl P^{(1)} + P^{(0)} \Big) \ti{\Phi},
\label{LP-VI2}
\end{align}
\end{subequations}
where the matrices $Q^{(1)}, P^{(1)}, Q^{(0)}, P^{(0)}$ are given by
\begin{subequations} \label{QP6}
\begin{align}
Q^{(1)}&=i\mathrm{diag}\,\big(-(\tau-1)\kappa_2-\tau\kappa_3,(\tau-1)\kappa_1-\kappa_3,\tau\kappa_1+\kappa_2\big),
 \label{Q6-1}\\
P^{(1)}&=i\mathrm{diag}\,\big(-\kappa_2-\kappa_3,\kappa_1,\kappa_1\big),
\label{P6-1}
\end{align}
and
\begin{align}
Q^{(0)} &= i
\begin{pmatrix}
-\gt_{23} & (\tau-1)^{-i\gr_{3}}v_{3}^{*} & -\tau^{i\gr_{2}}v_{2} \\
(\tau-1)^{i\gr_{3}}v_{3} & -\gt_{31} & v_{1}^{*} \\
-\tau^{-i\gr_{2}}v_{2}^{*} & v_{1} & -\gt_{12}
\end{pmatrix}, \label{Q6-0} \\
P^{(0)} &= i
\begin{pmatrix}
0 & (\tau-1)^{-1-i\gr_{3}}v_{3}^{*} & -\tau^{-1+i\gr_{2}}v_{2} \\
(\tau-1)^{-1+i\gr_{3}}v_{3} & 0 & 0 \\
-\tau^{-1-i\gr_{2}}v_{2}^{*} & 0 & 0
\end{pmatrix}. \label{P6-0}
\end{align}
\end{subequations}
In order to integrate the reduced system \eqref{RedSysVIa}, \eqref{RedSysVIb}
in terms of $P_6$ we compare the Fuchs--Garnier representation
\eqref{LP-VIgen} with the $3\times3$ Fuchs--Garnier representation
for $P_6$ obtained in \cite{H1994} and \cite{M2002}.
\begin{remark}
{\rm As noted earlier, the spectral parameter $k$ has been
introduced formally into Lax pair \eqref{3WRI-LP}. We made use of
this fact in extending the similarity reduction \eqref{SimRedVI} to
obtain a similarity reduction for the associated Lax pair. Here we
would like to illustrate that, although introduction of the
auxiliary spectral parameter $k$ is not absolutely necessary, in
this particular case it is an important ingredient of our
construction of the Fuchs--Garnier pair.

An alternate construction is also possible in which the spectral
variable $\gl$ is introduced without any dependence on $k$. Writing
$\ti{\Psi} = \Psi
\exp[ik(\kappa_{1}x_{1}+\kappa_{2}x_{2}+\kappa_{3}x_{3})]$ in
\eqref{3WRI-LP}, we then introduce the spectral variable as $\gl =
(x_{2} - x_{3})$ and follow the procedure described above. In fact,
we can just put $\kappa_1=\kappa_2=\kappa_3=0$ and $k=1$ in
\eqref{LP-VIgen}--\eqref{QP6}. As a result we arrive at the
following $3\times3$ matrix representation for the reduced system
\eqref{RedSysVIa}, \eqref{RedSysVIb}:
\begin{equation*}
\gl \frac{d\ti{\Phi}}{d\gl} = Q^{(0)} \ti{\Phi}, \quad
\frac{d\ti{\Phi}}{d\tau} = P^{(0)} \ti{\Phi},
\end{equation*}
where $Q^{(0)}, P^{(0)}$ are given in \eqref{QP6}. This is also a
$3\times3$ Fuchs--Garnier pair for the similarity solutions; one can
still get first integrals for system~\eqref{RedSysVIa}, \eqref{RedSysVIb}
as eigenvalues of $Q^{(0)}$, however all further information about the
solutions is hidden in a normalization of this system rather than
encoded in the monodromy structure. Therefore this system is
ineffective for further studying of the similarity solutions. }
\end{remark}
\subsection{Reduction of the $3\times3$ Fuchs--Garnier Pair to the $2\times2$ Pair in Jimbo-Miwa Form}
 \label{subsec:P6-FG-Laplace}

Now we simplify the notation and rewrite Fuchs--Garnier pair \eqref{LP-VIgen} in the following form:
\begin{subequations} \label{LP-P6}
\begin{align}
\Phi_{\gl} &= \Big( B^{6}_{1} + \frac{B^{6}_{0} - \mathrm{I}}{\gl} \Big) \Phi \label{LP6-1} \\
\Phi_{t} &= \Big( \gl M^{6}_{1} + M^{6}_{0} \Big) \Phi,
\label{LP6-2}
\end{align}
\end{subequations}
where the matrices $B^{6}_{1}, M^{6}_{1}, B^{6}_{0}, M^{6}_{0}$ are
given by
\begin{equation*}
B^{6}_{1} = \mathrm{diag}\, \big( t, 1, 0 \big), \quad M^{6}_{1} =
\mathrm{diag}\, \big( 1, 0, 0 \big),
\end{equation*}
and
\begin{equation*}
B^{6}_{0} =
\begin{pmatrix}
-\gt_{2} & \ti{w}_{3} & w_{2} \\
w_{3} & -\gt_{3} & \ti{w}_{1} \\
\ti{w}_{2} & w_{1} & -\gt_{1}
\end{pmatrix}, \quad
M^{6}_{0} =
\begin{pmatrix}
0 & (t-1)^{-1} \ti{w}_{3} & t^{-1} w_{2} \\
(t-1)^{-1} w_{3} & 0 & 0 \\
t^{-1} \ti{w}_{2} & 0 & 0
\end{pmatrix},
\end{equation*}
where $\{w_{j},\ti{w}_{j}\}$ are functions of $t$ and $\gt_{1},
\gt_{2}, \gt_{3}$ are arbitrary constants. Following \cite{M2002} we
assume that $0$ is one of eigenvalues of the matrix $B^{6}_{0}(t)$.
Note that this condition is a normalization of system~\eqref{LP-P6}
rather than a restriction. If we denote the eigenvalues of the
matrix $B^{6}_{0}(t)$ as $\mu_1$, $\mu_2$, and $\mu_3$, then we can
write:
\begin{equation*}
\mu_1=\frac12\Big(-\sum_{j=1}^3\gt_j+\gt_\infty\Big),\quad
\mu_2=\frac12\Big(-\sum_{j=1}^3\gt_j-\gt_\infty\Big),\quad\mu_3=0,
\end{equation*}
where $\gt_{\infty}$ is an arbitrary constant. We note that system
\eqref{LP6-1} coincides exactly with the system given in
\cite{M2002} if we make the gauge transformation $\Phi \mapsto
J\h{\Phi}$ where $J$ is the constant matrix
\begin{equation*}
J =
\begin{pmatrix}
0 & 1 & 0 \\
0 & 0 & 1 \\
1 & 0 & 0
\end{pmatrix}.
\end{equation*}

We omit writing the compatibility condition for pair \eqref{LP-P6},
which coincides with system \eqref{RedSysVIa}, \eqref{RedSysVIb}
rewritten in terms of variables $w_j,\tilde w_j$ (see equations~\eqref{Mazz-param2}),
because we do not use it. Instead, following \cite{M2002}, we briefly outline how this
pair can be mapped to the $2\times2$ Fuchs--Garnier pair for $P_6$
given by Jimbo and Miwa, which is defined by equation \eqref{P6-JMlam}.
As a result we obtain a parametrization of the 3WRI system in terms of
the solutions of $P_6$. We present this parametrization in the next
section.

We introduce the function $\ti{Y}(x,t)$ via the generalized Laplace
transform
\begin{equation} \label{P6-LT}
\Phi(\gl,t) = \int_{C} e^{\gl x} \ti{Y}(x,t) dx.
\end{equation}
Before substituting \eqref{P6-LT} into equations~\eqref{LP-P6}, it
is convenient to rewrite equation~\eqref{LP6-1} as follows,
\begin{equation*}
\gl \Phi_{\gl} = \big( \gl B^{6}_{1}(t) + B^{6}_{0}(t) - \mathrm{I}
\big) \Phi,
\end{equation*}
Assuming that the contour $C$ in \eqref{P6-LT} can be chosen to
eliminate any remainder terms that arise from integration-by-parts,
we find
\begin{equation*}
\big( B^{6}_{1}(t) - x \mathrm{I} \big) \frac{d\ti{Y}}{dx} =
B^{6}_{0}(t) \ti{Y},\qquad
\frac{d\ti{Y}}{dt}=-M_1^6\frac{d\ti{Y}}{dx}+M_0^6\ti{Y}.
\end{equation*}
Substituting the first equation obtained above into the second one
we obtain:
\begin{equation}
 \label{eq:FG-P6-ti-Y}
\frac{d\ti{Y}}{dx}=\big( B^6_1- xI\big)^{-1}B^6_0\ti{Y},\qquad
\frac{d\ti{Y}}{dt}=(-M_1^6\big( B^6_1-
xI\big)^{-1}B^6_0+M_0^6)\ti{Y},
\end{equation}
Since one of the eigenvalues of $B^6_0$, which are integrals of
motion, is 0, we can choose the Jordan form of $B_0^6$ such that all
the elements of its last column are zeroes. Denote such Jordan form
as $\h B_0^6$. If $B_0^6$ is diagonalizable, then $\h{B}^{6}_{0} =
\mathrm{diag}\, [\mu_{1}, \mu_{2}, 0 ]$. Define $G_0$, $\det G_0=1$,
as follows $ G_0^{-1}B_0^6(t_0)G_0=\h B_0^6 $ at some point
$t_0\neq0,1,\infty$. Now, define matrix $G$, as a solution of the
equation $ \frac{d}{dt}G=M_0^6G $ satisfying the initial data
$G(t_0)=G_0$. It is easy to observe that the compatibility
conditions for Fuchs--Garnier pair \eqref{LP-P6} implies that
$G^{-1}B_0^6(t)G=\h B_0^6$ holds for all $t$. We make the gauge
transformation $\ti Y=G\h Y$ in system~\eqref{eq:FG-P6-ti-Y}, to
find the following Fuchsian system for $\h Y$:
\begin{equation} \label{R3-Fuchsian}
\frac{d\h{Y}}{dx}=\left(\frac{\h{A}^6_0(t)}{x}+\frac{\h{A}^6_t(t)}{x-t}+\frac{\h{A}^6_1(t)}{x-1}\right)\h{Y},
\qquad\frac{d\h{Y}}{dt}=-\frac{\h{A}^6_t(t)}{x-t}\h{Y},
\end{equation}
where the $3\times3$ matrices $\h{A}^{6}_{j}$ all have the form
\begin{equation*}
\h{A}^{6}_{j} =
\begin{pmatrix}
* & * & 0 \\
* & * & 0 \\
* & * & 0
\end{pmatrix}.
\end{equation*}
Since the third column of each $\h{A}^{6}_{j}$ is zero, the system
for $\h{Y}$ reduces to a system for the first two components
\begin{equation} \label{P6-JM}
\frac{dY}{dx}=\left(\frac{A^6_0(t)}{x}+\frac{A^6_t(t)}{x-t}+\frac{A^6_1(t)}{x-1}\right)Y,\quad
\frac{dY}{dt}=-\frac{A^6_t(t)}{x-t}Y,\qquad Y =
\begin{pmatrix} \h{Y}_{1} \\ \h{Y}_{2} \end{pmatrix},
\end{equation}
and a quadrature for the third component. The eigenvalues of the
matrices $A^{6}_{0}$, $A^{6}_{t}$ and $A^{6}_{1}$ are $(\gt_{1},0)$,
$(\gt_{2},0)$ and $(\gt_{3},0)$, respectively, see \cite{M2002}.
Equation \eqref{P6-JM} is equivalent (up to gauge transformation) to
the $2\times2$ Fuchs--Garnier system for $P_6$ in the form given by
Jimbo and Miwa in \cite{JM1981II}. Now comparing with the Jimbo-Miwa
parametrization of the matrix elements of \eqref{P6-JM} by solutions
of $P_6$, we arrive at the parametrization for $w_j,\tilde w_j$
presented in the next section.

\subsection{Similarity Solution of 3WRI System in Terms of the Sixth Painlev\'{e} Equation}
 \label{subsec:P6-parametrization}

More details of the calculation explained in the previous section
can be found in \cite{M2002}. Here we present the final result, the
parametrization of the functions $w_j,\tilde w_j$ in terms of $P_6$
together with the corresponding reduction to get solutions of the
3WRI system.
\begin{subequations} \label{Mazz-param}
\begin{align}
w_{1}&=f\left(\frac{(t-1)y'-\theta_1(y-1)}{2y}+\frac{\gt_{3}(t-1)+(\gt_{\infty}-1)(y-1)}{2t}\right),\\
\ti{w}_{1}&=f^{-1}\left(-\frac{\gt_{3}y-ty'}{2(y-1)}+\frac{\gt_{1}t+(\gt_{\infty}-1)y}{2(t-1)}\right),\\
w_{2}&=\frac{g}{f}\left(-\frac{\gt_{2}y+ty'}{2(y-t)}-\frac{\gt_{1}+\gt_{\infty}y}{2(t-1)}+\frac{y(y-1)}
{2(t-1)(y-t)}\right),\\
\ti{w}_{2}&=\frac{f}{g}\left(\frac{t(t-1)y'-\theta_1(y-t)}{2y}-\frac{\gt_{2}(t-1)-
\gt_{\infty}(y-t)+y-1}{2}\right),\\
w_{3}&=g^{-1}\left(-\frac{\gt_{3}(y-t)+t(t-1)y'}{2(y-1)}+\frac{\gt_{2}t-\gt_{\infty}(y-t)+y}{2}\right),\\
\ti{w}_{3}&=g\left(-\frac{(t-1)y'+\gt_{2}(y-1)}{2(y-t)}+
\frac{\gt_{3}-\gt_{\infty}(y-1)}{2t}+\frac{y(y-1)}{2t(y-t)}\right),
\end{align}
\end{subequations}
where the functions $f=f(t)$ and $g=g(t)$ are the general solutions
of the following equations:
\begin{align}
 \label{eq:log-f}
\frac{d}{dt}\log{f}&=-\frac{y'}{2y(y-1)}-\frac{1+\gt_1-\gt_2+\gt_3}{2t(t-1)}+\frac{\gt_1}{2(t-1)y}+
\frac{\gt_3}{2t(y-1)},\\
\frac{d}{dt}\log{g}&=\frac{y'-1}{2(y-t)}-\frac{y'}{2(y-1)}+\frac{1-\gt_1+\gt_2-\gt_3}{2t}+\gt_2\left(\frac1{t-1}+
\frac1{2(y-t)}\right)\nonumber\\
 \label{eq:log-g}
&\hspace{10mm}+\gt_3\left(-\frac1{t(t-1)}+\frac1{2t(y-1)}\right),
\end{align}
and $y(t)$ is a solution of $P_{6}$ with
\begin{equation*}
\ga = \frac{(\gt_{\infty} - 1)^{2}}{2}, \quad \gb =
-\frac{\gt_{1}^{2}}{2}, \quad \gga = \frac{\gt_{3}^{2}}{2}, \quad
\gd = \frac{1 - \gt_{2}^{2}}{2},
\end{equation*}

In order to solve the reduced coupled 3WRI system \eqref{RedSysVIa}, \eqref{RedSysVIb}
in terms of $P_{6}$ we put $\kappa_{1} = 0, \kappa_{2} = 0, \kappa_{3} = i$
in \eqref{3WRI-LP}, and then compare matrix entries in the
Fuchs--Garnier pair \eqref{LP-VIgen} with $\tau = t$ with those in
system \eqref{LP-P6}. We obtain the following correspondence:
\begin{subequations} \label{Mazz-param2}
\begin{align}
iv_{1}(t) &= w_{1}(t),& -it^{i\gr_{2}}v_{2}(t) &= w_{2}(t),& i(t-1)^{i\gr_{3}}v_{3}(t) &= w_{3}(t),& \\
iv^{*}_{1}(t) &= \ti{w}_{1}(t),& -it^{-i\gr_{2}}v^{*}_{2}(t) &=
\ti{w}_{2}(t),& i(t-1)^{-i\gr_{3}}v^{*}_{3}(t) &= \ti{w}_{3}(t),&
\end{align}
and
\begin{equation}
i\gr_{1} = \gt_{1} - \gt_{3},\quad i\gr_{2} = \gt_{2} -
\gt_{1},\quad i\gr_{3} = \gt_{3} - \gt_{2}.
\end{equation}
\end{subequations}

Now we consider the physical reduction, i.e., assume that the star 
in
\eqref{RedSysVIa} denotes the complex conjugation. First of all we have
to impose the reduction on the formal monodromies:
\begin{equation*}
\gt_{1}=i\gt_{12},\qquad \gt_{2}=i\gt_{23},\qquad \gt_{3}=i\gt_{31},
\end{equation*}
where $\gt_{ik}\in\mathbb{R}$ and  $\gt_\infty\in\mathbb R$.
The solution $y(t)$ should be real for real $t$, and the functions $f$
and $g$ are as follows:
\begin{align*}
f(t)&=\frac{\sqrt{y}\,\sqrt{t}}{\sqrt{y-1}\sqrt{t-1}}
\left|\frac{t}{t-1}\right|^{\frac{\gt_{1}-\gt_{2}+\gt_{3}}{2}}
\exp{\left(\frac{\gt_{3}}2\int^t_{t_0}\frac{dt}{t(y-1)} +
\frac{\gt_{1}}2\int^t_{t_0}\frac{dt}{(t-1)y} + ic_1\right)},\\
g(t)&=\frac{\sqrt{y-t}\,\sqrt{t}}{\sqrt{y-1}}\,|t|^{\frac{\gt_{2}-\gt_{1}+\gt_{3}}2}|t-1|^{\gt_{2}-\gt_{3}}
\exp\left(\frac{\gt_{2}}2\int^t_{t_0}\frac{dt}{y-t} +
\frac{\gt_{3}}2\int^t_{t_0}\frac{dt}{t(y-1)}+ic_2\right),
\end{align*}
where the parameters $t_0,c_1,c_2\in\mathbb R$, the parameter $c_1\neq0$
only in the case if $\gt_{3}=\gt_{1}=0$, and $c_2\neq0$ if $\gt_{3}=\gt_{2}=0$.
Moreover, the solution of $P_6$ should satisfy the following condition: $0<t<1$ and $t<y(t)<1$.

\begin{remark}
{\rm We note that the parametrization that was adopted in
\cite{M2002} to write system \eqref{LP-P6} explicitly in terms of
$y$ where $y(t)$ is a solution of $P_{6}$ is not unique. Alternate
parameterizations have been identified by Boalch
\cite{B2005a}--\cite{B2005b} in his studies of $P_6$.
Since this system can be
mapped to the irregular $3\times3$ Lax pair of \cite{H1994} and
\cite{M2002} via the generalized Laplace transform, see equation
\eqref{R3-Fuchsian} above, it follows that these parameterizations
are equivalent to system \eqref{LP-P6} up to a gauge transformation. }
\end{remark}


\section{Similarity Reduction to the Fifth Painlev\'{e} Equation}
 \label{sec:P5}

We consider the following similarity reduction of the 3WRI system,
which was obtained in \cite{MW1989},
\begin{equation} \label{SimRedV}
\begin{split}
u_{1} = e^{-ix_{2}x_{3}} x_{3}^{\frac{i\gr}{2}} v_{1}, \quad u_{2} =
e^{ix_{3}x_{1}} x_{3}^{\frac{i\gr}{2}} v_{2}, \quad u_{3} = (x_{1} -
x_{2})^{-1+i\gr} v_{3},
\end{split}
\end{equation}
where $v_j=v_{j}(\tau)$ with
\begin{equation}
 \label{SimRedV-tau}
\quad \tau = (x_{1}-x_{2})x_{3},
\end{equation}
and $\gr$ is a real constant. Under these assumptions system
\eqref{3WRIcharac} reduces to the system of ODEs:
\begin{equation} \label{RedSysV}
\begin{split}
\tau^{1+i\gr} e^{i\tau} v_{1}' = iv_{2}^{*}v_{3}^{*}, \quad
\tau^{1+i\gr} e^{i\tau} v_{2}' = -iv_{3}^{*}v_{1}^{*}, \quad
\tau^{i\gr} e^{i\tau} v_{3}' = iv_{1}^{*}v_{2}^{*}.
\end{split}
\end{equation}
where prime denotes differentiation with respect to $\tau$. This
system was integrated in \cite{MW1989} in terms of an SD-function
and shown to be solvable in terms of the fifth Painlev\'{e}
equation~\eqref{P5}.

\begin{remark}{\rm
It is straightforward to generalize this similarity reduction to the
coupled case of the 3WRI system. One adds to \eqref{SimRedV} and
\eqref{RedSysV} the formally conjugated equations
\begin{equation*}
\begin{gathered}
u_1^*=e^{ix_2x_3}x_3^{-\frac{i\gr}2}v_1^*,\quad
u_2^*=e^{-ix_3x_1}x_3^{-\frac{i\gr}2}v_2^*,\quad
u_3^*=(x_1-x_2)^{-1-i\gr} v_3^*,\\
\tau^{1-i\gr} e^{-i\tau}{v_1^*}'=-iv_2v_3,\quad \tau^{1-i\gr}
e^{-i\tau}v_2*'=+iv_3v_1,\quad \tau^{-i\gr}
e^{-i\tau}v_3*'=-iv_1v_2.
\end{gathered}
\end{equation*}
Note that in the coupled case $\rho\in\mathbb C$ and the functions
$v_j$ and $v_j^*$ are not assumed to be complex conjugates. In the
most part of this Section we deal with the coupled 3WRI system and
turn back to the physical case at the end of
Subsection~\ref{subsec:P5-parametrization}.}
\end{remark}

\subsection{Fuchs--Garnier Pair for the Reduced System}

Following the approach outlined in the previous section we will use
\eqref{SimRedV} to construct a $3\times3$ Fuchs--Garnier pair for
the reduced system \eqref{RedSysV}. The pair is valid in the coupled
case also.

Consider Lax pair \eqref{3WRI-LP}. In this case we introduce the
spectral parameter in a different way comparing with the previous
section: instead of a scaled version of the spectral parameter $k$,
the new spectral parameter $\gl$ is defined in terms of the
dynamical variables, namely,
\begin{equation*}
\gl = (x_{1}+x_{2})x_{3}.
\end{equation*}
Setting $\kappa_{1} = \kappa_{2} = \kappa_{3} = 0$ in
\eqref{3WRI-LP} one proves that the solution of the Lax pair has the
following similarity form,
$$
\Psi=R(x_{j})\ti{\Phi}(\tau,\gl),\qquad R(x_{1},x_{2},x_{3}) =
\mathrm{diag}\, \Big( e^{ix_{1}x_{3}} x_{3}^{-i\gt_{23}},
e^{ix_{2}x_{3}} x_{3}^{-i\gt_{31}}, x_{3}^{-1-i\gt_{12}} \Big),
$$
where
\begin{equation} \label{theta-5}
\gt_{12} - \gt_{31}= -\frac{\gr}{2}, \quad \gt_{23} - \gt_{12} =
-\frac{\gr}{2}, \quad \gt_{31} - \gt_{23} = \gr.
\end{equation}
In terms of the new variables the Lax pair \eqref{3WRI-LP} becomes
\begin{equation*}
\begin{split}
\tau\ti{\Phi}_{\tau}+D_{1}\ti{\Phi}_{\gl}&=i\Big(-\tfrac{1}{2}(\gl-\tau)S_{2}+V_{1}\Big)\ti{\Phi}\\
\tau\ti{\Phi}_{\tau}+D_{2}\ti{\Phi}_{\gl}&=i\Big(-\tfrac{1}{2}(\gl+\tau)S_{1}+V_{2}\Big)\ti{\Phi},
\end{split}
\end{equation*}
where the matrices $D_{j}, S_{j}, V_{j}$ are given by
\begin{align*}
D_{1} &= \mathrm{diag}\, \big( -\tau, \gl, \tau \big),& D_{2} &= \mathrm{diag}\, \big( \gl, \tau, -\tau \big), \\
S_{1} &= \mathrm{diag}\, \big( 1, 0, 0 \big),& S_{2} &= \mathrm{diag}\, \big( 0, 1, 0 \big), \\
V_{1} &=
\begin{pmatrix}
0 & \tau^{-i\gr} e^{-i\tau} v_{3}^{*} & 0 \\
0 & \gt_{31} & -v_{1}^{*} \\
-\tau v_{2}^{*} & 0 & 0
\end{pmatrix},&
V_{2} &=
\begin{pmatrix}
\gt_{23} & 0 & v_{2} \\
\tau^{i\gr} e^{i\tau} v_{3} & 0 & 0 \\
0 & -\tau v_{1} & 0
\end{pmatrix}.
\end{align*}
After rearranging, the above system can be written as
\begin{subequations} \label{LP-Vgen}
\begin{align}
\ti{\Phi}_{\gl} &= \Big( \frac{Q^{(0)}}{\gl + \tau} + \frac{Q^{(1)}}
{\gl - \tau} + Q^{(2)} \Big) \ti{\Phi} \label{LP-V1} \\
\ti{\Phi}_{\tau} &= \Big( \frac{Q^{(0)}}{\gl + \tau} -
\frac{Q^{(1)}} {\gl - \tau} + P^{(2)} \Big) \ti{\Phi} \label{LP-V2}
\end{align}
\end{subequations}
where the matrices $Q^{(0)}, Q^{(1)}$ and $Q^{(2)}, P^{(2)}$ are
given by
\begin{equation*}
Q^{(0)} = i
\begin{pmatrix}
\gt_{23} & -\tau^{-i\gr} e^{-i\tau} v_{3}^{*} & v_{2} \\
0 & 0 & 0 \\
0 & 0 & 0
\end{pmatrix} \quad
Q^{(1)} = i
\begin{pmatrix}
0 & 0 & 0 \\
-\tau^{i\gr} e^{i\tau} v_{3} & \gt_{31} & -v_{1}^{*} \\
0 & 0 & 0
\end{pmatrix}
\end{equation*}
and
\begin{equation*}
Q^{(2)} = -\frac{i}{2}
\begin{pmatrix}
1 & 0 & 0 \\
0 & 1 & 0 \\
v_{2}^{*} & -v_{1} & 0
\end{pmatrix}, \quad
P^{(2)} = -\frac{i}{2}
\begin{pmatrix}
1 & -2\tau^{-1-i\gr} e^{-i\tau} v_{3}^{*} & 0 \\
-2\tau^{-1+i\gr} e^{i\tau} v_{3} & -1 & 0 \\
v_{2}^{*} & v_{1} & 0
\end{pmatrix}.
\end{equation*}

\subsection{Fuchs--Garnier Pairs for the Fifth Painlev\'{e} Equation}
 \label{subsec:P5-FG-Laplace}
The goal of this section is to establish a map between the
$3\times3$ Fuchs--Garnier pair \eqref{LP-Vgen} and the $2\times2$
Fuchs--Garnier pair for $P_5$ found by Jimbo and Miwa
\cite{JM1981II}. It is convenient to introduce the ``coupled
notation'' for matrix elements of the Fuchs--Garnier pair:
\begin{subequations}
 \label{eq:P5-FG-3x3}
\begin{align}
\Phi_{\gl}&=\left(\frac{B_1}{\gl+t}+\frac{B_2}{\gl-t}+\frac12I+B_3\right)\Phi,\label{LP5-1}\\
\Phi_t&=\left(\frac{B_1}{\gl+t}-\frac{B_2}{\gl-t}+M_\infty\right)\Phi,\label{LP5-2}
\end{align}
\end{subequations}
where
\begin{equation}
 \label{eq:P5-Bk-def}
B_1=
\begin{pmatrix}
\ti{m} & \ti{w_{3}} & w_{2} \\
0 & 0 & 0 \\
0 & 0 & 0
\end{pmatrix},\quad
B_2=
\begin{pmatrix}
0 & 0 & 0 \\
w_{3} & m & \ti{w_{1}} \\
0 & 0 & 0
\end{pmatrix},\quad
B_3=
\begin{pmatrix}
0 & 0 & 0 \\
0 & 0 & 0 \\
\ti{w_2}/2 & w_1/2 & -1/2
\end{pmatrix},
\end{equation}
where $I$ is the identity matrix and
\begin{equation}
 \label{eq:P5-Minfty-def}
M_\infty=
\begin{pmatrix}
1/2 & -t^{-1}\ti{w_{3}} & 0 \\
-t^{-1}w_{3} & -1/2 & 0 \\
\ti{w_2}/2 & -w_1/2 & 0
\end{pmatrix},
\end{equation}
where $\{w_{j},\ti{w_{j}}\}$ are all functions of $t$.

Compatibility of equations \eqref{LP5-1} and \eqref{LP5-2} gives the
following system of equations:
\begin{equation}
m' = 0, \quad \ti{m}' = 0,
\end{equation}
and
\begin{equation} \label{P5sys}
\begin{aligned}
tw_{1}' &= \ti{w_{2}}\ti{w_{3}},& t\ti{w_{1}}' &= -w_{2}w_{3}, \\
tw_{2}' &= -\ti{w_{1}}\ti{w_{3}},& t\ti{w_{2}}' &= w_{1}w_{3}, \\
tw_{3}' &= -[t - (m - \ti{m})] w_{3} - t \ti{w_{1}}\ti{w_{2}}, &
t\ti{w_{3}}' &= [t - (m - \ti{m})] \ti{w_{3}} + t w_{1}w_{2},
\end{aligned}
\end{equation}
where the primes denote derivatives by $t$.

As in the previous section, we are going to apply to the
Fuchs--Garnier pair the generalized Laplace transform. For this
purpose we rewrite equations ~\eqref{eq:P5-FG-3x3} in the
appropriate form with the coefficients linearly depending on the
spectral parameter, namely:
\begin{subequations}
 \label{eq:FG-P5-linear}
\begin{align}
 \label{eq:FG-P5-lambda-linear}
(\gl J_0+tJ)\Phi_\gl&=\big(\frac12(\gl J_0+tJ)+B\big)\Phi,\\
\label{eq:FG-P5-t-linear} (\gl J_0+tJ)\Phi_t&=\big((\gl
J_0+tJ)M+JB\big)\Phi,
\end{align}
\end{subequations}
where
\begin{subequations}
\begin{gather}
 \label{eq:P5-J0J}
J_0=
\begin{pmatrix}
1 & 0 & 0 \\
0 & 1 & 0 \\
0 & 0 & 0
\end{pmatrix}, \qquad
J=
\begin{pmatrix}
1 & 0 & 0 \\
0 & -1 & 0 \\
0 & 0 & 1
\end{pmatrix},\\
 \label{eq:P5-B}
B=B_1+B_2+tB_3=
\begin{pmatrix}
\ti{m} & \ti{w_{3}} & w_{2} \\
w_{3} & m & \ti{w_{1}} \\
\frac{t}{2}\ti{w_{2}} & \frac{t}{2}w_{1} & -\frac{t}{2}
\end{pmatrix},\\
\label{eq:P5-M} M=M_\infty-B_3=
\begin{pmatrix}
1/2 &- \ti{w_3}/t & 0 \\
-w_3/t & -1/2 & 0 \\
0 & -w_1 & 1/2
\end{pmatrix},
\end{gather}
\end{subequations}
We define the generalized Laplace transform as follows,
\begin{equation} \label{P5-LT}
\Phi(\gl,t)=\int_Ce^{\gl x/2}\ti{Y}(x,t)dx.
\end{equation}
Substituting it into equations \eqref{eq:FG-P5-lambda-linear} and
\eqref{eq:FG-P5-t-linear}, and assuming that the contour $C$ can be
suitably chosen to eliminate any remainder terms that arise from
integration-by-parts, we find
\begin{gather}
 \label{eq:P5-FG-3lambda}
(x-1)J_0\frac{d\ti{Y}}{dx}=\Big(\frac{t}2(x-1)J-(J_0+B)\Big)\ti{Y},\\
 \label{eq:P5-FG-3t}
(J_0+B)\frac{d\tilde
Y}{dt}=\Big(\frac{x-1}2(J+JB+tJM-tJ_0MJ)+J_0M(J_0+B)-\frac{d}{dt}B\Big)\tilde
Y,
\end{gather}
where in the derivation of equation \eqref{eq:P5-FG-3t} we used
equation~\eqref{eq:P5-FG-3lambda}. The third row of equation
\eqref{eq:P5-FG-3lambda} reads:
\begin{equation}
x \ti{Y}_{3} = \ti{w}_{2} \ti{Y}_{1} + w_{1} \ti{Y}_{2}.
\end{equation}
Using this relation to eliminate $\ti{Y}_{3}$ from
\eqref{eq:P5-FG-3lambda} and \eqref{eq:P5-FG-3t} we obtain the
following $2\times2$ system:
\begin{align}
 \label{eq:P5-FG-2lambda}
\frac{dY}{dx}& = \left( \frac{t}{2}
\begin{pmatrix}
1 & 0 \\
0 & -1
\end{pmatrix}
+ \frac{1}{x}
\begin{pmatrix}
w_{2}\ti{w}_{2} & w_{1}w_{2} \\
\ti{w}_{1}\ti{w}_{2} & w_{1}\ti{w}_{1}
\end{pmatrix}
- \frac{1}{x-1}
\begin{pmatrix}
w_{2}\ti{w}_{2} + \ti{m}+1 & \ti{w_{3}} + w_{1}w_{2}\\
w_{3} + \ti{w}_{1}\ti{w_{2}}& w_{1}\ti{w}_{1} + m +1
\end{pmatrix}\right) Y,\\
 \label{eq:P5-FG-2t}
\frac{dY}{dt}& = \left( \frac{x}{2}
\begin{pmatrix}
1 & 0 \\
0 & -1
\end{pmatrix}
+ \frac{1}{t}
\begin{pmatrix}
0 & -\ti{w_3}\\
-w_{3}& 0
\end{pmatrix}
\right) Y,\qquad Y = \begin{pmatrix} \ti{Y}_{1} \\ \ti{Y}_{2}
\end{pmatrix}.
\end{align}
The Fuchs--Garnier pair \eqref{eq:P5-FG-2lambda},
\eqref{eq:P5-FG-2t} coincides (up to a simple gauge transformation)
with the Fuchs--Garnier pair for $P_5$ by Jimbo--Miwa (see
\cite{JM1981II}, equations (C.38), (C.39)).

\subsection{Parametrization of Solutions in Terms of $P_{5}$}
 \label{subsec:P5-parametrization}

In order to parameterize the general solution of
system~\eqref{P5sys} by the (general) solution of $P_5$ we compare
the parametrization $2\times2$ Fuchs--Garnier representations for
this system obtained above \eqref{eq:P5-FG-2lambda},
\eqref{eq:P5-FG-2t} with the one by Jimbo and Miwa \cite{JM1981II}.

First of all we notice that system~\eqref{P5sys} admits first
integrals
\begin{equation} \label{integrals5a}
\begin{split}
&m={\rm const},\quad\ti{m}-m =\theta_\infty, \quad w_{1}\ti{w_{1}} + w_{2}\ti{w_{2}} = \theta_0 \\
&w_1w_2w_3+\ti{w_1}\ti{w_2}\ti{w_3}+w_3\ti{w_3}+\frac{\theta_\infty}2(w_2\ti{w_2}-w_1\ti{w_1})=
\frac{\theta_1^2-\theta_0^2-\theta_\infty^2}4,
\end{split}
\end{equation}
where $\gt_{0},\gt_{1},\gt_{\infty}$ are arbitrary constants, which
have a sense of formal monodromies of the normalized solution $Y$ of
system~\eqref{eq:P5-FG-2lambda}, \eqref{eq:P5-FG-2t}. ``Normalized''
here means that we make a gauge transformation of $Y$ which puts all
matrices in \eqref{eq:P5-FG-2lambda} into the traceless form. The
notation of the formal monodromies coincides with those from the
Jimbo-Miwa work \cite{JM1981II}. The first integral $m$ cannot be
expressed via the monodromies because the normalized version of
equation \eqref{eq:P5-FG-2lambda} depends only on the difference
$\tilde m-m$, rather than on $m$ and $\tilde m$ separately.

Motivated by the parametrization used by Jimbo--Miwa we define the
functions
\begin{equation} \label{P5-eqn}
\begin{split}
y(t)&=\frac{w_2\ti{w}_2(\ti{w}_3+w_1w_2)}{(w_2\ti{w}_2-(\gt_0+\gt_1-\gt_\infty)/2)w_1w_2},\quad
z(t)=w_2\ti{w_2}-\gt_0,\quad u(t)=-\frac{w_1}{\tilde w_2}.
\end{split}
\end{equation}
It follows from system \eqref{P5sys} and parametrization
\eqref{integrals5a} that $y$, $z$, $u$ satisfy the following system
of nonlinear ODEs:
\begin{subequations} \label{C.40}
\begin{align}
ty' &= ty - 2z(y-1)^{2} - \left( \frac{\gt_{0} - \gt_{1} +
\gt_{\infty}}{2} \right) (y-1)^{2} +
(\gt_{0} + \gt_{1}) (y-1), \label{y-dash} \\
 \label{eq:P5-z}
tz' &= yz \left(z + \frac{\gt_{0} - \gt_{1} + \gt_{\infty}}{2}
\right) - \frac{1}{y} (z + \gt_{0}) \left(z +
\frac{\gt_{0} + \gt_{1} + \gt_{\infty}}{2} \right), \\
 \label{eq:P5-u}
t(\log
u)'&=-2z-\theta_0+y\left(z+\frac{\theta_0-\theta_1+\theta_\infty}2\right)+
\frac1y\left(z+\frac{\theta_0+\theta_1+\theta_\infty}2\right).
\end{align}
\end{subequations}
System \eqref{C.40} coincides with the system (C.40) in
\cite{JM1981II}. Eliminating $z$ from the first equation and
substituting it into the second one we find that $y(t)$ satisfies
the general $P_{5}$ equation \eqref{P5} with the coefficients:
\begin{equation}
 \label{eq:P5-coefficients-thetas}
\ga=\frac12\left(\frac{\gt_0-\gt_1+\gt_\infty}2\right)^2,\quad
\gb=-\frac12\left(\frac{\gt_0-\gt_1-\gt_\infty}2\right)^{2},\quad
\gga=1-\gt_0-\gt_1,\quad\gd = -1/2.
\end{equation}

Now we find the converse formulae, namely, the functions
$\{w_{j}(t),\ti{w}_{j}(t)\}$ in terms of $y(t)$ and $z(t)$ and
$u(t)$. Using \eqref{P5-eqn} we obtain the following representations
for $\{w_{j}(t),\ti{w}_{j}(t)\}$:
\begin{equation} \label{P5-wfunctions}
\begin{aligned}
&\begin{aligned} w_{1} &= -f z(z + \gt_{0}),&
w_{2} &= \frac{1}{gz},&\\
\ti{w_{1}} &= \frac{1}{f(z+\gt_{0})},& \ti{w}_{2} &= g
z(z+\gt_{0}),&
\end{aligned} \\
&\begin{aligned}
w_{3} &= \frac{g}{f} \left( \frac{1}{y} \Big( z + \frac{\gt_{0} + \gt_{1} + \gt_{\infty}}{2} \Big) - z \right), \\
\ti{w_{3}}&= -\frac{f}{g} \left( y \Big( z + \frac{\gt_{0} - \gt_{1}
+ \gt_{\infty}}{2} \Big) - (z + \gt_{0}) \right),
\end{aligned}
\end{aligned}
\end{equation}
where instead of one function $u(t)$ we are forced to introduce two
functions $f(t)$ and $g(t)$, such that $u(t)=f(t)/g(t)$. The
additional function appears as a result of an extra ``gauge
freedom'' in the $3\times3$ Fuchs--Garnier pair compared to the
$2\times2$ one. System \eqref{P5sys} implies, that the functions
$f(t)$ and $g(t)$ satisfy the following equations:
\begin{equation}
 \label{eq:P5-fg-yz}
\begin{aligned}
t(\log{f})' &= -\frac{tz'}{z + \gt_{0}} - \frac{tz'}{2z} +
\frac{1}{2}
\left( y \Big( z + \frac{\gt_{0} - \gt_{1} + \gt_{\infty}}{2} \Big) - (z + \gt_{0}) \right)\\
&\hspace{20mm} + \frac{z + \gt_{0}}{2z}
\left( \frac{1}{y} \Big( z + \frac{\gt_{0} + \gt_{1} + \gt_{\infty}}{2} \Big) - z \right), \\
t(\log{g})' &= -\frac{tz'}{z} - \frac{tz'}{2(z + \gt_{0})} -
\frac{z}{2(z + \gt_{0})}
\left( y \Big( z + \frac{\gt_{0} - \gt_{1} + \gt_{\infty}}{2} \Big) - (z + \gt_{0}) \right) \\
&\hspace{20mm}-\frac12\left(\frac1{y}\Big(z+\frac{\gt_0+\gt_1+
\gt_\infty}2\Big)-z\right).
\end{aligned}
\end{equation}
These expressions can be simplified by introducing the function
$\gs(t)$, following the work of Jimbo and Miwa in \cite{JM1981II}.
In our notation, we define $\gs(t)$ as
\begin{equation} \label{P5-sigma}
\gs(t) = w_{3}\ti{w}_{3} + tw_{1}\ti{w}_{1} + \frac{(\gt_{0} +
\gt_{\infty})^{2} - \gt_{1}^{2}}{4}.
\end{equation}
Then, using the fourth identity in \eqref{integrals5a}, we find
\begin{subequations} \label{P5-fgfunctions}
\begin{align}
t(\log{f})' &= -\frac{tz'}{z + \gt_{0}} - \frac{tz'}{2z} + \frac{1}{2z} \Big( \gs + (t + \gt_{\infty})\gs' \Big), \\
t(\log{g})' &= -\frac{tz'}{z} - \frac{tz'}{2(z + \gt_{0})} -
\frac{1}{2(z + \gt_{0})} \Big( \gs + (t + \gt_{\infty})\gs' \Big).
\end{align}
\end{subequations}
We note that the function $\gs(t)$ satisfies the following two
important equations:
\begin{equation}
 \label{eq:P5-sigma-prime}
 \frac{d\gs}{dt}=-z(t),
\end{equation}
which can be proved by the differentiation of equation
\eqref{P5-sigma}, and
\begin{multline} \label{P5-sigmaODE}
t^2\left(\frac{d^2\gs}{dt^2}\right)^2=\left(\gs-(\theta_\infty+2\theta_0+t)\frac{d\gs}{dt}+
2\Big(\frac{d\gs}{dt}\Big)^2\right)^2\\
-4\frac{d\gs}{dt}\left(\frac{d\gs}{dt}-\theta_0\right)
\left(\frac{d\gs}{dt}-\frac{\theta_0+\theta_1+\theta_\infty}2\right)
\left(\frac{d\gs}{dt}-\frac{\theta_0-\theta_1+\theta_\infty}2\right).
\end{multline}
Equation \eqref{P5-sigmaODE} can be verified in the following way.
Substituting the parametrization \eqref{P5-wfunctions} into equation
\eqref{P5-sigma}, we express the $\sigma$-function in terms of $y$
and $z$. We then couple the resulting expression with
equation~\eqref{eq:P5-z} and use \eqref{eq:P5-sigma-prime} to
eliminate $z$. Then, summing up and subtracting these equations one
finds the two equations:
\begin{equation}
 \label{eq:y-1/y-sigma}
\begin{aligned}
t\sigma''+(\sigma-(\theta_\infty+2\theta_0+t)\sigma'+2\sigma'^2)&=2\frac{\sigma'-\theta_0}y
\left(\sigma'-\frac{\theta_0+\theta_1+\theta_\infty}2\right),\\
t\sigma''-(\sigma-(\theta_\infty+2\theta_0+t)\sigma'+2\sigma'^2)&
=-2\sigma'y
\left(\sigma'-\frac{\theta_0-\theta_1+\theta_\infty}2\right),
\end{aligned}
\end{equation}
The compatibility condition of these equations is equivalent to
\eqref{P5-sigmaODE}. Thus $\sigma(t)$ is the so-called
$SD$-function.

We are now ready to solve the reduced 3WRI system \eqref{RedSysV} in
terms of $P_{5}$. The Fuchs--Garnier pairs \eqref{LP-Vgen} and
\eqref{eq:P5-FG-3x3} are related by the change of variables $\gl
\mapsto i\gl, \tau \mapsto it$. By comparing matrix entries between
\eqref{LP-Vgen} and \eqref{eq:P5-FG-3x3} we get the following
correspondence
\begin{subequations}
 \label{eq:P5-vj-wj}
\begin{align}
v_{1}(\tau) &= -w_{1}(t),& v_{2}(\tau) &= -iw_{2}(t),& v_{3}(\tau) &= i\tau^{-i\gr}e^{-i\tau} w_{3}(t),& \\
v^{*}_{1}(\tau) &= i\ti{w}_{1}(t),& v^{*}_{2}(\tau) &=
\ti{w}_{2}(t),& v^{*}_{3}(\tau) & = i\tau^{i\gr}e^{i\tau}
\ti{w}_{3}(t),&
\end{align}
\begin{equation}
i\gt_{23} =\tilde m,\quad i\gt_{31} =m, \quad\gr =
\theta_{31}-\theta_{23}=i(\tilde m-m)=i\gt_{\infty}.
\end{equation}
\end{subequations}
where the functions $\{ w_{j}, \ti{w}_{j} \}$ are given in terms of
$y$ and $z$ by equations \eqref{P5-wfunctions} and
\eqref{eq:P5-fg-yz}. These formulae define the general similarity
solution for the coupled case of the 3WRI system. We note that as
follows from equations \eqref{eq:P5-sigma-prime} and
\eqref{eq:y-1/y-sigma} the solution of the coupled system can be
presented in terms one function $\sigma(t)$.

To give the general similarity solution in the physical case of the
3WRI system, i.e.~where we prove that the functions $v_j(\tau)$ and
$v_j^*(\tau)$, $j = 1,2,3$, are complex conjugates for real $\tau$,
it is necessary to present the solution solely in terms of the
function $\gs(t)$. To achieve this we impose the following
conditions on the parameters
\begin{equation}
 \label{eq:P5-real}
t,\theta_0,\theta_1,\theta_\infty\in i\mathbb R.
\end{equation}
Introducing the notation $\tilde\sigma(\tau)=\sigma(t)$, where
$\tau=it$, we note that $\ti{\gs}$ satisfies an ODE analogous to
\eqref{P5-sigmaODE} which, by condition \eqref{eq:P5-real}, will
have real coefficients. It follows that we can take the general
solution $\ti{\gs}(\tau)$ of this equation to be real. After making
the change of variables $t=-i\tau$, $\theta_{\infty}=-i\rho$,
$\theta_0=-i\rho_0$, $\theta_1=-i\rho_1$ with
$\tau,\rho,\rho_0,\rho_1 \in \mathbb{R}$, we define the functions
$z(t)$ and $y(t)$ by equations \eqref{eq:P5-sigma-prime} and any one
of \eqref{eq:y-1/y-sigma}, respectively. Then, we use
\eqref{P5-fgfunctions} to obtain the following expressions for the
functions $f$ and $g$:
\begin{equation}
\begin{aligned}
&\frac1{f(z+\gt_0)}=z^{1/2}\exp{\left(-i\int_{\tau_0}^\tau\frac{\ti{\gs}-(\tau+\gr)\ti{\gs}'}{2\ti{\gs}'}d\tau+
if_0\right)},\\
&\frac1{gz}=(z+\gt_0)^{1/2}\exp{\left(i\int_{\tau_0}^\tau\frac{\ti{\gs}-(\tau+\gr)}{2(\ti{\gs}'+\gr_0)}d\tau+
ig_0\right)},
\end{aligned}
\end{equation}
Notice that the function $z(t)=-i\tilde\sigma'(\tau)$ is a pure
imaginary function of $\tau$. Assume that $e^{i\pi/2}z>0$ and
$e^{-i\pi/2}(z+\theta_0)>0$.
It is straightforward now to observe
from equations \eqref{P5-wfunctions} and \eqref{eq:P5-vj-wj} that
the functions $v_j$ and $v_j^*$ for $j=1,2$ are indeed conjugates
under conditions \eqref{eq:P5-real}. To see that the same is true
for the functions $v_3$ and $v_3^*$, one has to employ additionally
equations~\eqref{eq:y-1/y-sigma}.

\subsection{Alternate Reduction to the $2\times2$ Fuchs--Garnier Pair for $P_{5}$}
 \label{subsec:P5alt}

In this section we present an alternate reduction of the $3\times3$
system \eqref{eq:P5-FG-3x3} to the Jimbo--Miwa version of the
Fuchs--Garnier system for $P_5$ \cite{JM1981II}, by making use of
suitable gauge transformations rather than the generalized Laplace
transform. The key observation is that the parameter $m$ in matrix
$B$ (see equation~\eqref{eq:P5-B}) can be chosen such that its
determinant vanishes for all values of $t$. Indeed,
\begin{equation*}
\det\,B=w_1w_2w_3+\ti{w_1}\ti{w_2}\ti{w_3}+w_3\ti{w_3}-\ti{m}w_1\ti{w}_1-mw_2\ti{w_2}-m\ti{m},
\end{equation*}
is the first integral of system~\eqref{P5sys} by virtue of
equations~\eqref{integrals5a}. Moreover, by putting
\begin{equation} \label{eq:P5alt-m-def}
m = -\frac{\gt_{0} + \gt_{1} + \gt_{\infty}}{2},
\end{equation}
one finds that $\det\,B$ coincides with the difference of the
l.-h.s. and r.-h.s. of the last integral in \eqref{integrals5a} and
thus vanishes for all $t$\footnote{Parameter $\theta_1$ is defined
in \eqref{integrals5a} up to the sign, therefore one can change
$\theta_1\to-\theta_1$ in the definition of m
in~\eqref{eq:P5alt-m-def}.}. It follows that $B$ has eigenvalues
$(\mu_{1}, \mu_{2}, 0)$, where $\mu_k=\mu_k(t)$, $k=1,2$. Moreover,
on the {\bf general} solutions\footnote{\label{foot1}We assume that
the general situation holds in this section.} of
system~\eqref{P5sys} all eigenvalues are pairwise different, thus
there exists an invertible matrix $G=G(t)$ such that
\begin{equation}
 \label{eq:P5alt-tildeB}
G^{-1} B G=
\begin{pmatrix}
\mu_{1} & 0 & 0 \\
0 & \mu_{2} & 0 \\
0 & 0 & 0
\end{pmatrix}
\equiv\tilde B.
\end{equation}
In Section~\ref{subsec:Okamoto-BT} we give explicit expressions for
the eigenvalues $\mu_{1}(t), \mu_{2}(t)$ and the diagonalizing
matrix $G(t)$ in terms of the functions $\{ w_{j}, \ti{w}_{j} \}$,
$j = 1,2,3$.

We now make the gauge transformation $\Phi = G \tilde{\mathcal{Y}}$
in system \eqref{eq:FG-P5-lambda-linear}, \eqref{eq:FG-P5-t-linear}
to obtain
\begin{equation}
 \label{eq:P5alt-3x3linear}
\frac{d\tilde{\mathcal{Y}}}{d\gl}= \left(\frac12I+\tilde
A_3+\frac{\tilde A_2}{\gl - t}+\frac{\tilde A_1}{\gl +
t}\right)\tilde{\mathcal{Y}},\qquad \frac{d\tilde{\mathcal{Y}}}{dt}
= \left(\tilde M_\infty-\frac{\tilde A_2}{\gl - t}+\frac{\tilde
A_1}{\gl + t} \right)\tilde{\mathcal{Y}},
\end{equation}
where
\begin{align}
 \label{eq:P5alt-Ak}
\tilde A_k&=G^{-1}B_kG=G^{-1}I_kBG=G^{-1}I_kG\tilde B,\\
\tilde M_\infty&=G^{-1}\left(M_\infty G-\frac{d}{dt}G\right),
 \label{eq:P5alt-tildeMinfty}
\end{align}
the matrices $B_k$, $M_\infty$, and $\tilde B$ are defined by
equations~\eqref{eq:P5-Bk-def}, \eqref{eq:P5-Minfty-def}, and
\eqref{eq:P5alt-tildeB}, respectively, and, for each $k=1,2,3$, the
matrix $I_k$ has only one nonzero element, which is the $k$-th
element on the diagonal, more precisely,
$I_k\equiv\{\delta_{ik}\delta_{kj}\}_{i,j=1}^{i,j=3}$, where
$\delta_{nm}$ is Kroneker's delta.

Note that from the first equation~\eqref{eq:P5alt-Ak} follows that
$\rm{rank}\,\tilde A_k=1$. Let us prove that the first two elements
in the third columns of the matrices $\tilde M_\infty$ and $\tilde
A_k$ are zeroes:
\begin{equation}
 \label{eq:P5alt-13-23-zeroes}
\tilde M_\infty[1,3]=\tilde M_\infty[2,3]=\tilde A_k[1,3]=\tilde
A_k[2,3]=A_k[3,3]=0,\qquad k=1,2,3.
\end{equation}
For the matrix elements of $\tilde A_k$
equations~\eqref{eq:P5alt-13-23-zeroes} are an immediate consequence
of the second equation~\eqref{eq:P5alt-Ak}. For the matrix $\tilde
M_\infty$ it follows from the compatibility condition of
system~\eqref{eq:P5alt-3x3linear}, $\tilde A_k'=[\tilde
M_\infty,A_k]$ , where the brackets denote the matrix commutator,
the corresponding structure of the matrices $\tilde A_k$, see
equations~\eqref{eq:P5alt-13-23-zeroes}, and the fact that ${\rm
rank}\,B=2$ according to the assumption in footnote~\ref{foot1}.
\begin{equation}
 \label{eq:P5alt-MinftyAk-block}
\tilde M_\infty=\left(\!\!\!
\begin{array}{cc}
\begin{array}{c}\fbox{\Large$\hat M_\infty$}\vspace{2pt}\\
*\quad*\end{array}\hspace{-4pt}&\hspace{-10pt}\begin{array}{c}0\\0\\\phantom{'}*\phantom{'}\end{array}\\
\end{array}
\!\!\!\right),\qquad \tilde A_k=\left(\!\!\!
\begin{array}{cc}
\begin{array}{c}\fbox{\Large$\hat A_k$}\vspace{2pt}\\
*\quad*\end{array}\hspace{-4pt}&\hspace{-10pt}\begin{array}{c}0\\0\\0\end{array}\\
\end{array}
\!\!\!\right),\qquad k=1,2,3.
\end{equation}
Note that $\tilde M_\infty[3,3]\neq0$.

Thus the structure of the matrices $\tilde M_\infty$ and $\tilde
A_k$ \eqref{eq:P5alt-MinftyAk-block} implies that
system~\eqref{eq:P5alt-3x3linear} can be reduced to a system in
$2\times2$ matrices for the first two components of
$\tilde{\mathcal{Y}}$, while the third component can be found in
terms of the first two via a quadrature. The reduced $2\times2$
system looks the same as system~\eqref{eq:P5alt-3x3linear} the only
difference is that tildes are changed by hats.

To proceed let us notice that the eigenvalues of the matrices
$\tilde A_1$, $\tilde A_2$, and $\tilde A_3$, equal $\{\tilde m,
0,0\}$, $\{m,0,0\}$, and $\{-1/2,0,0\}$, correspondingly; it follows
from the definition of $\tilde A_k$ and $B_k$ see the first equation
in \eqref{eq:P5alt-Ak} and equations \eqref{eq:P5-Bk-def},
respectively. The structure of matrices $\tilde A_k$
in~\eqref{eq:P5alt-MinftyAk-block} implies that the matrices $\hat
A_1$, $\hat A_2$, and $\hat A_3$, have the following eigenvalues:
$\{\tilde m, 0\}$, $\{m,0\}$, and $\{-1/2,0\}$, respectively.
Therefore, the matrix $1/2I+\hat A_3$ has eigenvalues $\{1/2,0\}$
and there exists an invertible matrix $H$, such that,
\begin{equation}
 \label{eq:P5alt-H-def}
H^{-1} (\tfrac{1}{2}I + \h{A}_{2}) H =
\begin{pmatrix}
\tfrac{1}{2} & 0 \\
0 & 0
\end{pmatrix},
\end{equation}
Thus, the function
$$
{\mathcal Y}\equiv H^{-1}\hat{\mathcal Y}=H^{-1}
\begin{pmatrix}
\tilde{\mathcal Y}_{11}& \tilde{\mathcal Y}_{12} \\
\tilde{\mathcal Y}_{21} & \tilde{\mathcal Y}_{22}
\end{pmatrix},
$$
where $\hat{\mathcal Y}$ is the corresponding main minor of any
solution of system~\eqref{eq:P5alt-3x3linear}, solves the
Fuchs--Garnier pair:
\begin{equation}
 \label{eq:P5alt-JM1}
\frac{d{\mathcal Y}}{d\gl}=\left(
\begin{pmatrix}
\tfrac{1}{2} & 0 \\
0 & 0
\end{pmatrix}+
\frac{{\mathcal A}_2}{\gl-t}+\frac{{\mathcal
A}_1}{\gl+t}\right){\mathcal Y},\qquad \frac{d{\mathcal
Y}}{dt}=\left( {\mathcal D}_\infty- \frac{{\mathcal
A}_2}{\gl-t}+\frac{{\mathcal A}_1}{\gl+t}\right){\mathcal Y},
\end{equation}
where
\begin{equation}
 \label{eq:P5alt-Dinfinity}
{\mathcal A}_k=H^{-1}\hat
A_kH\quad\text{for}\;\;k=1,2\quad\text{and}\quad {\mathcal
D}_\infty=H^{-1}\big(\hat M_\infty H-\frac{d}{dt}H\big).
\end{equation}
We remark that ${\mathcal D}_\infty$ is a diagonal matrix for any
$H$ in \eqref{eq:P5alt-H-def}. Indeed, the compatibility condition
for system~\eqref{eq:P5alt-JM1} implies,
$$
\big[{\mathcal D}_\infty,
\begin{pmatrix}
\tfrac{1}{2} & 0 \\
0 & 0
\end{pmatrix}\big]=0.
$$
Finally, to put the Fuchs--Garnier pair \eqref{eq:P5alt-JM1}, into
the Jimbo-Miwa form for $P_5$ \cite{JM1981II}, we have to make the
following transformation:
$$
{\mathcal
Y}(\lambda,t)=e^{\lambda/4}x^{\ti{m}/2}(x-1)^{m/2}D(t)Y(x,t),\qquad\lambda=2tx-t,
$$
where $x$ is the new spectral parameter and $D$ is a diagonal matrix
depending only on $t$ defined as follows,
\begin{equation}
 \label{eq:P5alt-D}
D^{-1}\frac{d}{dt}D={\mathcal D}_\infty+\frac1t{\rm
diag}\,({\mathcal A}_1+{\mathcal A}_2)-\frac14
\begin{pmatrix}
1&0\\
0&-1
\end{pmatrix}.
\end{equation}
The function $Y$ solves the system
\begin{equation}
 \label{eq:P5alt-JM-true}
\frac{dY}{dx}=\left( \frac{t}2\sigma_3
+\frac{A_0}x+\frac{A_1}{x-1}\right)Y,\qquad \frac{dY}{dt}=\left(
\frac{x}2\sigma_3 +\frac1t{\rm offdiag}(A_0+A_1) \right)Y,
\end{equation}
where the notation ${\rm offdiag}(\cdot)$ means the off-diagonal
part of the corresponding matrix, i.e., the matrix where the
diagonal elements are substituted by zeroes,
$$
\sigma_3=
\begin{pmatrix}
1&0\\
0&-1
\end{pmatrix},\qquad
A_0=D^{-1}{\mathcal A}_1D-\frac{\tilde m}2I,\qquad
A_1=D^{-1}{\mathcal A}_2D-\frac{m}2I.
$$
The matrices $A_k$ obey the following relations:
$$
{\rm tr}\,A_0={\rm tr}\,A_1=0,\quad \det\,A_0=-\frac{\tilde
m^2}4,\quad \det\,A_1=-\frac{m^2}4,\qquad {\rm
diag}(A_0+A_1)=-\frac{\theta_1-\theta_0}2\sigma_3.
$$
The results for the traces and determinants of $A_0$ and $A_1$ can
be deduced from the corresponding results for the matrices
${\mathcal A_1}$ and ${\mathcal A_2}$:
$$
{\rm tr}\,{\mathcal A}_1=\tilde m,\quad{\rm tr}\,{\mathcal
A}_2=m,\quad \det\,{\mathcal A}_1=\det\,{\mathcal A}_2=0,
$$
which are proved above. To prove the formula for the diagonal part
of $A_0+A_1$, we note that actually the following more general
formula is valid,
\begin{equation}
 \label{eq:P5alt-A0+A1}
A_0+A_1=D^{-1}H^{-1}
\begin{pmatrix}
\mu_1&0\\
0&\mu_2
\end{pmatrix}
HD+
\begin{pmatrix}
0&0\\
0&t/2
\end{pmatrix}
-\frac{m+\tilde m}2I.
\end{equation}
To prove identity~\eqref{eq:P5alt-A0+A1} one has to start with the formula
$B_1+B_2=B-B_3$ and follow the construction presented in this
section. Then use formula~\eqref{eq:P5Okamoto-H} for $H$ given in
Subsection~\ref{subsec:Okamoto-BT}, to prove that the diagonal part
of equation~\eqref{eq:P5alt-A0+A1} equals $(\theta_0+\frac{m+\tilde
m}2)\sigma_3$. Finally recall the choice of $m$ in
\eqref{eq:P5alt-m-def} and the equation for $\tilde m$ in
\eqref{integrals5a}.

By the way, since the trace of l.-h.s. of \eqref{eq:P5alt-A0+A1}
equals $0$, we find that $\mu_1+\mu_2=m+\tilde m-t/2$, which is
consistent with equation~\eqref{eq:P5Okamoto-mu-quadratic} of
Subsection~\ref{subsec:Okamoto-BT}.

\subsection{An Okamoto-type B\"{a}cklund Transformation for $P_{5}$}
 \label{subsec:Okamoto-BT}
In Subsections~\ref{subsec:P5-FG-Laplace} and \ref{subsec:P5alt} we
found two different reductions of the $3\times3$ Fuchs--Garnier
pair~\eqref{eq:P5-FG-3x3} to the $2\times2$ Fuchs--Garnier pair
of the Jimbo-Miwa type, namely, \eqref{eq:P5-FG-2lambda},
\eqref{eq:P5-FG-2t} and \eqref{eq:P5alt-JM-true}. In this Subsection
we present some details of the calculations related with the
reduction scheme of the previous Subsection. Using them the
interested reader can follow the same scheme as in
Subsection~\ref{subsec:P5-parametrization} to get an alternate
parametrization of the similarity reduction \eqref{SimRedV},
\eqref{SimRedV-tau} of the 3WRI system in terms of solutions of
$P_5$. We, however, proceed in a different way: we find a specific
Okamoto-type B\"acklund transformation for $P_5$ (see
Appendix~\ref{app:A} equation~\eqref{aeq:O}) together with the
generating integral transformation for solutions of the
Fuchs--Garnier pair.

So we begin with the explicit formulae for the objects introduced in
the previous Subsection: The diagonalizing matrix $G(t)$ in
\eqref{eq:P5alt-tildeB} is taken as
$$
\mathit{G}=\left({\begin{array}{ccc}
m\,{w_{2}} - {\tilde w_{1}}\,{\tilde w_{3}} - {w_{2}}\,{\mu _{1}}
&m\,w_2-\tilde w_1\,\tilde w_3-w_2\,\mu_2&m\,w_2-\tilde w_1\,\tilde w_3\\
\mathit{\tilde m}\,\tilde w_1-w_2\,w_3-\tilde w_1\,\mu_1&\mathit{\tilde m}\,\tilde w_1-
{w_{2}}\,{w_{3}} - {\tilde w_{1}}\,{\mu _{2}}
&\mathit{\tilde m}\,{\tilde w_{1}} - {w_{2}}\,{w_{3}} \\
 - {\displaystyle \frac {t}{2}} ({\mu _{2}} + {\theta _{0}}+ {\displaystyle \frac {t}{2}} )
 & - {\displaystyle \frac {t}{2}} ({\mu _{1}} + {\theta _{0}} + {\displaystyle \frac {t}{2} } )
 & - ({\mu _{1}} + {\displaystyle \frac {t}{2}} )\,({\mu _{2} } + {\displaystyle\frac {t}{2}} ) -
 {\displaystyle \frac {t}{2} } {\theta _{0}}
\end{array}}
 \right)
$$
where $\mu_1$ and $\mu_2$ are solutions of the following quadratic
equation,
\begin{equation}
\label{eq:P5Okamoto-mu-quadratic}
\mu^2-\big(m+\ti{m}-\frac{t}2\big)\mu-\big(w_3\ti{w}_3+\frac{t}2(m+\ti{m}+\gt_0)-m\ti{m}\big)=0.
\end{equation}
In the general situation all three eigenvalues $\{\mu_1,\mu_2,0\}$
are different,
$$
\det\,G=((m-\tilde m)w_2\tilde w_1+w_2^2w_3-\tilde w_1^2\tilde
w_3)(\mu_2-\mu_1)\mu_1\mu_2\neq0.
$$
The diagonalizing matrix $H(t)$ in \eqref{eq:P5alt-H-def} is taken
as
\begin{equation}
 \label{eq:P5Okamoto-H}
H=\left(
\begin{array}{cc}
1&-{\displaystyle\frac{\mu _2}{\mu_1}}\\
-{\displaystyle\frac{(2\,\mu _2+t+2\,\theta _0)\,\mu _1}
{(2\,\mu_1+t+2\,\theta _0)\,\mu_2}} &1
\end{array}
 \right),\qquad
\det\,H={\displaystyle\frac{2\,(\mu _1-\mu _2)}{2\,\mu_1+t+2\,\theta
_0}}.
\end{equation}
An important auxiliary object is the diagonal matrix $D$, the
logarithmic derivative of which is defined in
equation~\eqref{eq:P5alt-D}. A nontrivial ingredient of the formula
in \eqref{eq:P5alt-D} is the diagonal matrix $\mathcal{D}_\infty$
defined in the second equation in \eqref{eq:P5alt-Dinfinity}. Using
MAPLE code and following the algorithm of
Subsection~\ref{subsec:P5alt} one confirms that $\mathcal{D}_\infty$
is, indeed, the diagonal matrix. This calculation at the same time
gives extremely complicated expressions for the diagonal elements:
$\mathcal{D}_\infty[1,1]$, $\mathcal{D}_\infty[2,2]$. We were not
able to find a concise expression for them. At the same time it is
not complicated to find an expression for the logarithmic derivative
of the ratio $D_{11}/D_{22}$ of the diagonal elements of $D$, or,
equivalently, the difference
$\mathcal{D}_\infty[1,1]-\mathcal{D}_\infty[2,2]$, see below.

Using the formulae for $G(t)$, $H(t)$ and $D(t)$ given above we
obtain the following expressions for the matrices $A_{0}$, $A_{1}$
in \eqref{eq:P5alt-JM-true}:
\begin{align*}
A_1=\left(\begin{array}{c} {\displaystyle\frac{(w_2\,w_3-\tilde
m\,\tilde w_1)\,((\theta_0+\tilde m)\,w_2+\tilde w_1\,\tilde w_3)}
{(m-\tilde m)\,w_2\,\tilde w_1+w_2^2\,w_3-\tilde w_1^2\,\tilde w_3}}-
{\displaystyle\frac{\tilde m}2}\vspace{4pt}\\
{\displaystyle\frac{D_{11}\,\mu_1\,(\tilde m\,(\theta_0+m)\,\tilde
w_1 -\theta_0\,w_2\,w_3-\tilde w_1\,w_3\,\tilde w_3)\,
((\theta_0+\tilde m)\,w_2+\tilde w_1\,\tilde w_3)}
{D_{22}\,\mu_2\,(\mu_2+\theta_1)\, ((m-\tilde m)\,w_2\,\tilde
w_1+w_2^2\,w_3-\tilde w_1^2\,\tilde w_3)}}
\end{array} \right.\\
\left.\begin{array}{c}
{\displaystyle\frac{D_{22}\,\mu_2\,(\mu_2+\theta_1)\,w_2\,(w_2\,w_3-\tilde
m\,\tilde w_1)}
{D_{11}\,\mu_1\,((m-\tilde m)\,w_2\,\tilde w_1+w_2^2\,w_3-\tilde w_1^2\,\tilde w_3)}}\vspace{4pt}\\
{\displaystyle\frac {w_2(\tilde m\,(\theta_0+m)\,\tilde w_1
-\theta_0\,w_2\,w_3-\tilde w_1\,w_3\,\tilde w_3)} {(m-\tilde
m)\,w_2\,\tilde w_1+w_2^2\,w_3-\tilde w_1^2\,\tilde
w_3}}-{\displaystyle\frac{\tilde m}2}
\end{array}\right),
\end{align*}
\begin{align*}
A_2=\left(\begin{array}{c} {\displaystyle\frac{(m\,w_2-\tilde
w_1\,\tilde w_3)\, ((\theta_0+m)\,\tilde w_1+w_2\,w_3)}
{(m-\tilde m)\,w_2\,\tilde w_1+w_2^2\,w_3-\tilde w_1^2\,\tilde w_3}}-
{\displaystyle\frac{m}2}\vspace{4pt}\\
-{\displaystyle\frac{D_{11}\,\mu_1\,(m\,(\theta_0+\tilde
m)\,w_2-\theta_0\,\tilde w_1\,\tilde w_3- w_2\,w_3\,\tilde w_3)\,
((\theta_0+m)\,\tilde w_1+w_2\,w_3)}
{D_{22}\,\mu_2\,(\mu_2+\theta_1)\,((m-\tilde m)\,w_2\,\tilde
w_1+w_2^2\,w_3-\tilde w_1^2\,\tilde w_3)}}
\end{array} \right.\\
\left.\begin{array}{c}
{\displaystyle\frac{D_{22}\,\mu_2\,(\mu_2+\theta_1)\,\tilde
w_1\,(m\,w_2-\tilde w_1\,\tilde w_3)}
{D_{11}\,\mu_1\,((m-\tilde m)\,w_2\,\tilde w_1+w_2^2\,w_3-\tilde w_1^2\,\tilde w_3)}}\vspace{4pt}\\
-{\displaystyle\frac{\tilde w_1\,(m\,(\theta_0+\tilde
m)\,w_2-\theta_0\,\tilde w_1\,\tilde w_3-w_2\,w_3\,\tilde w_3)}
{(m-\tilde m)\,w_2\,\tilde w_1+w_2^2\,w_3-\tilde w_1^2\,\tilde
w_3}}-{\displaystyle\frac{m}2}
\end{array}\right),
\end{align*}
\begin{equation}
 \label{eq:P5Okamoto-A1+A2}
A_1+A_2=\left(\begin{array}{cc}
\displaystyle{\frac{\theta_0-\theta_1}2}&
-\displaystyle{\frac{D_{22}\,\mu_2\,(2\,\mu_1+2\,\theta_0+t)}
{2\,D_{11}\,\mu_1}}\vspace{4pt}\\
\displaystyle{\frac{D_{11}\,\mu_1\, (2\,\mu_2+2\,\theta_0+t)}
{2\,D_{22}\,\mu_2}}& -\displaystyle{\frac{\theta_0-\theta_1}2}
\end{array}\right).
\end{equation}
Equation~\eqref{eq:P5Okamoto-A1+A2} is obtained as a sum of the
matrices $A_0$ and $A_1$ presented above. However, we used
identities for $\mu_1$ and $\mu_2$ following from
equation~\eqref{eq:P5Okamoto-mu-quadratic} to simplify the off
diagonal elements. The same formula can be obtained in a different
way: from equation~\eqref{eq:P5alt-A0+A1} with matrix $H$ in
\eqref{eq:P5Okamoto-H}.

Now we compare the Jimbo-Miwa parametrization of
system~\eqref{eq:P5alt-JM-true} (\cite{JM1981II} equation
(C.38))\footnote{For the convenience of the reader this
parametrization is presented in equation
\eqref{aeq:JM-P5parametrization} in Appendix~\ref{app:A}.} with the
one obtained above. To differentiate from the solution of $P_5$ that
already appeared in Subsection~\ref{subsec:P5-parametrization} we
adopt the hat notation: the solution of $P_5$, $\hat y=\hat y(t)$,
associated functions $\hat z=\hat z (t)$ and $\hat u=\hat u(t)$, see
system~\eqref{C.40}, and the corresponding monodromies
$\hat\theta_0$, $\hat\theta_1$, $\hat\theta_\infty$, which we obtain
in this section. Thus we arrive at the following equations for the
formal monodromies:
\begin{equation}
 \label{eq:P5Okamoto-hat-thetas}
\hat\theta_0=-\tilde m=\frac{\theta_0+\theta_1-\theta_\infty}2,\quad
\hat\theta_1=m=-\frac{\theta_0+\theta_1+\theta_\infty}2,\quad
\hat\theta_\infty=\theta_1-\theta_0,
\end{equation}
The choice of the signs for $\hat\theta_0$ and $\hat\theta_1$ in
equations~\eqref{eq:P5Okamoto-hat-thetas} are in our hands (see
Appendix~\ref{app:A}). After we fixed the signs we obtain equations
for the $P_5$ functions:
\begin{equation}
 \label{eq:P5Okamoto-hat-z}
\begin{aligned}
\hat z&=-\frac{(\tilde m\tilde w_1-w_2w_3)((\theta_0+\tilde
m)w_2+\tilde w_1\tilde w_3)}
{(m-\tilde m)\tilde w_1w_2+w_2^2w_3-\tilde w_1^2\tilde w_3},\\
\hat z+\hat\theta_0&=-\frac{w_2(\tilde m(\theta_0+m)\tilde
w_1-\theta_0w_2w_3-\tilde w_1w_3\tilde w_3)}
{(m-\tilde m)\tilde w_1w_2+w_2^2w_3-\tilde w_1^2\tilde w_3},\\
\hat z+\frac{\hat\theta_0+\hat\theta_1+\hat\theta_\infty}2&=
-\frac{\tilde w_1(m(\theta_0+\tilde m)w_2-\theta_0\tilde w_1\tilde
w_3-w_2w_3\tilde w_3)}
{(m-\tilde m)\tilde w_1w_2+w_2^2w_3-\tilde w_1^2\tilde w_3},\\
\hat z+\frac{\hat\theta_0-\hat\theta_1+\hat\theta_\infty}2&=
-\frac{(mw_2-\tilde w_1\tilde w_3)((\theta_0+m)\tilde w_1+w_2w_3)}
{(m-\tilde m)\tilde w_1w_2+w_2^2w_3-\tilde w_1^2\tilde w_3},\\
\end{aligned}
\end{equation}
\begin{align}
 \label{eq:P5Okamoto-hat-y}
\hat y&=\frac{\tilde w_1(\tilde m(\theta_0+m)\tilde
w_1-\theta_0w_2w_3-\tilde w_1w_3\tilde w_3)}
{(\tilde m\tilde w_1-w_2w_3)((\theta_0+m)\tilde w_1+w_2w_3)},\\
 \label{eq:P5Okamoto-hat-u}
\hat u&=-\frac{D_{22}\mu_2(\mu_2+\theta_1)}{D_{11}\mu_1}
\frac{(\tilde m\tilde w_1-w_2w_3)}{(\tilde m(\theta_0+m)\tilde
w_1-\theta_0w_2w_3-\tilde w_1w_3\tilde w_3)}.
\end{align}
The formulae \eqref{eq:P5Okamoto-hat-z} arise from the comparison of
the different matrix elements, of course, all of them are
equivalent.

We can use now the methodology of
Subsection~\ref{subsec:P5-parametrization} to invert equations
\eqref{eq:P5Okamoto-hat-z}--\eqref{eq:P5Okamoto-hat-u} to get a
parametrization of the similarity solutions of 3WRI system in terms
of ``hat'' $P_5$ functions. However, there is much more sense to
rewrite these equations in terms of the ``uncovered'' $P_5$
functions obtained in Subsection~\ref{subsec:P5-parametrization} by
exploiting equations~\eqref{P5-wfunctions}. In this way we obtain
the Okamoto-type B\"acklund transformation:
\begin{equation}
 \label{eq:P5Okamoto-transformations-yzu}
\hat z=z+\frac{\theta_0-\theta_1+\theta_\infty}2,\qquad \hat
y=\displaystyle{\frac{yz}{z+\frac{\theta_0+\theta_1+\theta_\infty}2}},\qquad
\hat
u=\frac{D_{22}\mu_2(\mu_2+\theta_1)}{D_{11}\mu_1(z(y-1)-\theta_0)}.
\end{equation}
The first two formulae here represent the Okamoto transformation for
$P_5$. We complete them by the reference that $y$ solves \eqref{P5}
for the coefficients \eqref{eq:P5-coefficients-thetas}, while the
function $\hat y$ is the solution of \eqref{P5} for the following
set of the coefficients:
\begin{equation}
 \label{eq:P5Okamoto-coefficients-natural}
\begin{aligned}
\hat\alpha&=\frac12\left(\frac{\hat\theta_0-\hat\theta_1+\hat\theta_\infty}2\right)^2=\frac{\theta_1^2}2,&
\hat\beta&=-\frac12\left(\frac{\hat\theta_0-\hat\theta_1-\hat\theta_\infty}2\right)^2=-\frac{\theta_0^2}2,\\
\hat\gamma&=1-\hat\theta_0-\hat\theta_1=1+\theta_\infty,&
\delta&=-\frac12.
\end{aligned}
\end{equation}

Let us now consider the function $\hat u$. First of all, notice that
we can use the formula for the logarithmic derivative of $u$ in
terms of $y$ and $z$ \eqref{eq:P5-u} to get the corresponding
transformation for the logarithmic derivative of $\hat u$,
\begin{equation}
 \label{eq:P5-Okamoto-lderivative-hat-u}
t\frac{d}{dt}\log\,\hat
u=t\frac{d}{dt}\log\,u-\frac{\theta_0-\theta_1+\theta_\infty}2(y+1)
\end{equation}
The functions $\mu_1$ and $\mu_2$ are solutions of the quadratic
equation:
\begin{gather*}
\mu^2+\Big(\theta_0+\theta_1+\frac t2\Big)\mu
+\frac{t\theta_1}2+\frac{(\theta_0+\theta_1)^2-\theta_\infty^2}4-\\
\frac1y\left(\left(z+\frac{\theta_0-\theta_1+\theta_\infty}2\right)
(y-1)-\frac{\theta_0+\theta_1-\theta_\infty}2\right)
\left(z(y-1)-\frac{\theta_0+\theta_1+\theta_\infty}2\right)=0.
\end{gather*}
Although we do not present explicit formulae for the matrix $D$,
the logarithmic derivative of $D_{11}/D_{22}$ or the difference
$\mathcal{D}_\infty[1,1]-\mathcal{D}_\infty[2,2]$ can be found explicitly
by using the third equation in \eqref{eq:P5Okamoto-transformations-yzu}
and equation \eqref{eq:P5-Okamoto-lderivative-hat-u}.

\begin{remark}\label{rem:P5Okamoto-parametrization}{\rm
We note that in case we take the Fuchs--Garnier pair in Jimbo-Miwa
parametrization and substitute the the functions $y$, $z$, and $u$,
by $\hat y$, $\hat z$, and $\hat u$, then we get a more natural
parametrization of the Fuchs--Garnier pair with the $P_5$
functions: in this parametrization each formal monodromy is
responsible for the corresponding coefficient of $P_5$ cf.
\eqref{eq:P5Okamoto-coefficients-natural}. Note that Jimbo-Miwa
parameterizations for all other Painlev\'e equations \cite{JM1981II}
is similar to the one we are proposing in this remark, so in a sense
we are proposing the ``true''
Jimbo-Miwa parametrization for $P_5$.
Explicitly this parametrization is presented at the end of
Appendix~\ref{app:A}}
\end{remark}

Of course, the method we use here to obtain the Okamoto
transformation allows us to get the corresponding transformation for
the solutions of the Fuchs-Garnier pairs. We denote as $\hat Y(x,t)$
the solution of the Fuchs--Garnier pair \eqref{eq:P5alt-JM-true} of
Subsection~\ref{subsec:P5alt} so that to make this notation
consistent with the notation for $P_5$ functions introduced above.
The function $Y(x,t)$, in the following formula, is the solution of
the Fuchs--Garnier pair \eqref{eq:P5-FG-2lambda},
\eqref{eq:P5-FG-2t}\footnote{To put system \eqref{eq:P5-FG-2lambda},
\eqref{eq:P5-FG-2t} into the standard traceless form one has to make
an additional scalar gauge 
transformation, $Y\to
x^{\theta_0/2}(x-1)^{1+\theta_1/2}Y$.} of
Subsection~\ref{subsec:P5-FG-Laplace}. The formula relating these
functions reads,
\begin{equation}
 \label{eq:P5Okamoto-integral}
\hat Y(x,t)=x^{-\frac{\tilde
m}2}(x-1)^{-\frac{m}2}D(t)^{-1}\int\limits_Ce^{t(x-\frac12)(\tilde
x-\frac12)} \left(P+\frac1{\tilde x}Q\right)Y(\tilde x,t)\,d\tilde
x,
\end{equation}
where the numbers $\tilde m$ and $m$ are defined in
\eqref{eq:P5Okamoto-hat-thetas}, matrix $D(t)$ - in
\eqref{eq:P5alt-D}, contour $C$ is the same as in
Subsection~\ref{subsec:P5-FG-Laplace}, and
\begin{gather*}
P=H^{-1}F,\qquad Q=H^{-1}\begin{pmatrix}
f_{13}\\
f_{23}
\end{pmatrix}\cdot
\begin{pmatrix}
\frac2tB_{31}&\frac2tB_{32}
\end{pmatrix},\\
\text{where $F$ is the $2\times2$ submatrix of}\qquad
G^{-1}=\left(\!\!\!
\begin{array}{cc}
\begin{array}{c}\fbox{\Large $F$}\vspace{2pt}\\
\vspace{-5pt}*\quad*\end{array}\hspace{-4pt}&\hspace{-10pt}
\begin{array}{c}f_{13}\\f_{23}\\\phantom{'}*\phantom{'}\end{array}\\
\end{array}
\!\!\!\right),
\end{gather*}
$G^{-1}$ is the inverse of $G$, defined in the beginning of this
subsection, and $B_{3k}$ $k=1,2$ are the matrix elements of $B$
\eqref{eq:P5-B}. The explicit expressions for $P$ and $Q$ are as
follows:
\begin{align*}
P&=D_\mu\begin{pmatrix} -\displaystyle{\frac{\tilde m\tilde
w_1-w_2w_3} {(m-\tilde m)\tilde w_1w_2+w_2^2-\tilde w_1^2\tilde
w_3}}& \displaystyle{\frac{mw_2-\tilde w_1\tilde w_3}
{(m-\tilde m)\tilde w_1w_2+w_2^2-\tilde w_1^2\tilde w_3}}\vspace{4pt}\\
\displaystyle{\frac{\tilde m(\theta_0+m)\tilde
w_1-\theta_0w_2w_3-\tilde w_1w_3\tilde w_3} {(m-\tilde m)\tilde
w_1w_2+w_2^2-\tilde w_1^2\tilde w_3}}&
-\displaystyle{\frac{m(\theta_0+\tilde m)w_2-\theta_0\tilde
w_1\tilde w_3-w_2w_3\tilde w_3} {(m-\tilde m)\tilde
w_1w_2+w_2^2-\tilde w_1^2\tilde w_3}}
\end{pmatrix}\!,\\
Q&=D_\mu\begin{pmatrix}
0&0\\
\tilde w_2&w_1
\end{pmatrix},
\qquad D_\mu=\frac1{\mu_1-\mu_2}\begin{pmatrix}
\displaystyle{\frac{\mu_2+\theta_1}{\mu_1}}&0\\
0&\displaystyle{\frac1{\mu_2}}
\end{pmatrix}.
\end{align*}


\section{Similarity Reduction to the Fourth Painlev\'{e} Equation}
 \label{sec:P4}

The following similarity reduction of 3WRI system~\eqref{3WRIcharac}
was found in \cite{MW1989}:
\begin{equation} \label{SimRedIV}
u_j = e^{i\phi_{j}}v_j(\tau),\qquad j=1,2,3,\qquad \tau = x_{1} +
x_{2} + x_{3},
\end{equation}
where
\begin{equation} \label{SimRedIVb}
\begin{split}
\phi_{1} &= \gr x_{3} + \tfrac{1}{2} x_{3}^{2} + 2x_{2}x_{3} +
\tfrac{1}{2} \gr^{2}, \quad
\phi_{2} = \gr x_{3} + \tfrac{1}{2} x_{3}^{2} + 2x_{3}x_{1} + \tfrac{1}{2} \gr^{2}, \\
\phi_{3} &= 2\gr (x_{1} + x_{2}) + (x_{1} + x_{2})^{2}, \\
\end{split}
\end{equation}
and $\gr$ is a real constant. Under these conditions system
\eqref{3WRIcharac} reduces to the system of ODEs:
\begin{equation} \label{RedSysIV}
\begin{split}
e^{i\phi} v_{1}' = iv_{2}^{*}v_{3}^{*}, \quad e^{i\phi} v_{2}' =
iv_{3}^{*}v_{1}^{*}, \quad e^{i\phi} v_{3}' = iv_{1}^{*}v_{2}^{*},
\end{split}
\end{equation}
where
\begin{equation}
 \label{eq:P4-phi}
\phi = \phi_{1} + \phi_{2} + \phi_{3} = (\tau + \gr)^{2},
\end{equation}
and prime denotes differentiation with respect to $\tau$. This
system was integrated in \cite{MW1989} in terms of SD-functions by
splitting real and imaginary parts of the equations. These
SD-functions were shown to be related with the fourth Painlev\'{e}
functions~\eqref{P4}.

\begin{remark}{\rm
As usual, it is straightforward to generalize this similarity
reduction to the coupled case of the 3WRI system. One adds to
\eqref{SimRedIV} and \eqref{RedSysIV} the formally conjugated
equations:
\begin{equation*}
 \label{eq:P4-coupled}
u_j^*=e^{-i\phi_j}v_j^*,\qquad e^{-i\phi} {v_j^*}'=-iv_kv_l,
\end{equation*}
respectively, where $(j,k,l)$ is any cyclic permutation of
$(1,2,3)$, $\phi_j$ and $\phi$ are defined in \eqref{SimRedIVb} and
\eqref{eq:P4-phi}, respectively, with $\rho\in\mathbb C$. As usual
in the coupled case the functions $v_j$ and $v_j^*$ are not assumed
to be complex conjugates. In the most part of this Section we deal
with the coupled 3WRI system and turn back to the physical case at
the end of Subsection~\ref{subsec:P4-parametrization}.}
\end{remark}

\subsection{A $3\times3$ Fuchs--Garnier Pair for the Reduced System}

Following the approach of the previous sections, we will use
\eqref{SimRedIV} to construct a Fuchs--Garnier pair which is valid
for both the coupled and physical cases of the reduced system
\eqref{RedSysIV}.

Consider the Lax pair \eqref{3WRI-LP}. Instead of the spectral
parameter $k$ we define the spectral parameter $\gl$ in the
following way
\begin{equation} \label{lambdaIV}
\gl = x_{1} - x_{2}.
\end{equation}
Since the spectral parameter is already defined we put
$\kappa_1=\kappa_2=\kappa_3=0$, and by the direct substitution prove
that $\Psi(x_{j},k) = R(x_{j})\ti{\Phi}(\tau,\gl)$, where $R(x_{j})$
is given by
\begin{equation*}
R(x_{1},x_{2},x_{3}) = \mathrm{diag}\, \Big( e^{i\phi_{2}},
e^{-i\phi_{1}}, 1 \Big).
\end{equation*}
In the new variables the Lax pair takes the form:
\begin{equation} \label{Phi4}
\begin{split}
\ti{\Phi}_{\tau} + D_{1} \ti{\Phi}_{\gl} &= i\big( -(\gl - \tau) S_{2} + V_{1} \big) \ti{\Phi} \\
\ti{\Phi}_{\tau} + D_{2} \ti{\Phi}_{\gl} &= i\big( -(\gl + \tau)
S_{1} + V_{2} \big) \ti{\Phi},
\end{split}
\end{equation}
where the matrices $D_{j}, S_{j}, V_{j}$ are given by
\begin{align*}
D_{1} &= \mathrm{diag}\, \big( -1, 0, 1 \big),& D_{2} &= \mathrm{diag}\, \big( 0, 1, -1 \big), \\
S_{1} &= \mathrm{diag}\, \big( 1, 0, 0 \big),& S_{2} &= \mathrm{diag}\, \big( 0, 1, 0 \big), \\
V_{1} &=
\begin{pmatrix}
0 & -e^{-i\phi} v_{3}^{*} & 0 \\
0 & \tfrac{1}{2}\gr & -v_{1}^{*} \\
-v_{2}^{*} & 0 & 0
\end{pmatrix},&
V_{2} &=
\begin{pmatrix}
-\tfrac{1}{2}\gr & 0 & v_{2} \\
e^{i\phi} v_{3} & 0 & 0 \\
0 & v_{1} & 0
\end{pmatrix}.
\end{align*}
After rearranging, the above system can be written in the form
\begin{subequations} \label{LP-IVgen}
\begin{align}
\ti{\Phi}_{\gl} &= \Big( \gl Q^{(1)} + Q^{(0)} \Big) \ti{\Phi} \label{LP-IV1} \\
\ti{\Phi}_{\tau} &= \Big( \gl P^{(1)} + P^{(0)} \Big) \ti{\Phi},
\label{LP-IV2}
\end{align}
\end{subequations}
where the matrices $Q^{(1)}, P^{(1)}, Q^{(0)}, P^{(0)}$ are given by
\begin{align*}
Q^{(1)} &= -i\mathrm{diag}\, \big( 1, -1, 0 \big),& P^{(1)} &=
-i\mathrm{diag}\, \big( 1, 1, 0 \big),
\end{align*}
and
\begin{equation*}
Q^{(0)} = i
\begin{pmatrix}
-(\tau + \gr) & e^{-i\phi} v_{3}^{*} & v_{2} \\
e^{i\phi} v_{3} & -(\tau + \gr) & v_{1}^{*} \\
-\tfrac{1}{2} v_{2}^{*} & -\tfrac{1}{2} v_{1} & 0
\end{pmatrix}, \quad
P^{(0)} =
\begin{pmatrix}
-(\tau + \gr) & 0 & v_{2} \\
0 & (\tau + \gr) & -v_{1}^{*} \\
-\tfrac{1}{2} v_{2}^{*} & \tfrac{1}{2} v_{1} & 0
\end{pmatrix}.
\end{equation*}

\subsection{Fuchs--Garnier Pairs for the Fourth Painlev\'{e} Equation}
 \label{subsec:P4-FG-Laplace}

We now consider the Fuchs--Garnier pair \eqref{LP-IVgen} in more
detail: we introduce variables $\{ w_{j}, \ti{w}_{j} \}$, $j=1,2,3$,
to emphasize the ``coupled character'' of the system under
consideration and write
\begin{subequations} \label{LP-P4}
\begin{align}
\Phi_{\gl} &= \Big( \gl B^{4}_{1} + B^{4}_{0} \Big) \Phi \label{LP4-1} \\
\Phi_{\tau} &= \Big( \gl M^{4}_{1} + M^{4}_{0} \Big) \Phi,
\label{LP4-2}
\end{align}
\end{subequations}
where the matrices $B^{4}_{1}, M^{4}_{1}, B^{4}_{0}, M^{4}_{0}$ are
given by
\begin{equation} \label{B4-1}
B^{4}_{1} = -i\mathrm{diag}\, \big( 1, -1, 0 \big), \quad M^{4}_{1}
= -i\mathrm{diag}\, \big( 1, 1, 0 \big),
\end{equation}
and
\begin{equation} \label{B4-0}
B^{4}_{0} =
\begin{pmatrix}
-i(\tau + \gr) & \ti{w_{3}} & w_{2} \\
w_{3} & -i(\tau + \gr) & \ti{w_{1}} \\
\ti{w_{2}} & w_{1} & 0
\end{pmatrix}, \quad
M^{4}_{0} =
\begin{pmatrix}
-i(\tau + \gr) & 0 & w_{2} \\
0 & i(\tau + \gr) & -\ti{w_{1}} \\
\ti{w_{2}} & -w_{1} & 0
\end{pmatrix},
\end{equation}
and $\{w_{j},\ti{w}_{j}\}$ are all functions of $\tau$.
Compatibility of equations \eqref{LP4-1} and \eqref{LP4-2} gives the
following system of equations
\begin{equation}
 \label{P4sys}
\begin{aligned}
w_{1}' &= \ti{w_{2}}\ti{w_{3}},& \ti{w_{1}}' &= -w_{2}w_{3},\\
w_{2}' &= \ti{w_{1}}\ti{w_{3}},& \ti{w_{2}}' &= -w_{1}w_{3}, \\
w_{3}'&=2i(\tau+\gr)w_{3}-2\ti{w_{1}}\ti{w_{2}},&\ti{w_{3}}'&=-2i(\tau+\gr)\ti{w_{3}}+2w_{1}w_{2}.
\end{aligned}
\end{equation}

Our next goal is to use the generalized Laplace transform
\eqref{P6-LT} to construct the map between the $3\times3$
Fuchs--Garnier pair \eqref{LP-P4} and the one in $2\times2$ matrices
for $P_4$ found by Jimbo and Miwa \cite{JM1981II}.

Substituting the formula for $\Phi(\gl,\tau)$ from equation
\eqref{P6-LT} into equations \eqref{LP4-1} and \eqref{LP4-2}, and
assuming that the contour $C$ is suitably chosen to eliminate any
remainder terms that arise from integration-by-parts, we find
\begin{align} \label{P4-LT1}
B^{4}_{1} \frac{d\ti{Y}}{dx}& = \big( -x I + B^{4}_{0} \big) \ti{Y},\\
 \label{P4-LT2}
\frac{d\tilde Y}{dt}&=\big(ixB^4_1+M^4_0-iB^4_1B^4_0\big)\tilde Y,
\end{align}
where for the derivation of equation~\eqref{P4-LT2} we used the
identity $i\big(B^4_1\big)^2=M^4_1$ (see equations~\eqref{B4-1}) and
the matrices $B^4_0, M^4_0$ are given in \eqref{B4-0}. We note that,
because the diagonal matrix $B^{4}_{1}$ has a zero in the (33)
entry, the third rows of these equations give the following
relationship between the components of $\ti{Y}$
\begin{equation*}
x \ti{Y}_{3} = \ti{w_{2}} \ti{Y}_{1} + w_{1} \ti{Y}_{2},\qquad
\frac{d}{dt}\tilde Y_3=\tilde w_2\tilde Y_1-w_1\tilde Y_2.
\end{equation*}
Using the first equation above to eliminate $\ti{Y}_{3}$ from
equations \eqref{P4-LT1} and \eqref{P4-LT2} we arrive at the
following $2\times2$ system
\begin{equation} \label{P4-JM}
\begin{aligned}
\frac{d\h{Y}}{dx}&=\left(xA^4_2+A^4_1(\tau)+\frac{A^4_0(\tau)}x\right)\h{Y},\\
\frac{d\hat Y}{dt}&=\left(xA^4_2+A^4_1(\tau)-i(\tau+\rho)A^4_2\right)\hat
Y, \quad \h{Y} = \begin{pmatrix} \ti{Y}_{1} \\ \ti{Y}_{2}
\end{pmatrix},
\end{aligned}
\end{equation}
where
\begin{equation*}
A^{4}_{0}(\tau) =
\begin{pmatrix}
-w_{2}\ti{w_{2}} & -w_{2}w_{1} \\
\ti{w_{2}}\ti{w_{1}} & w_{1}\ti{w_{1}}
\end{pmatrix}, \;
A^{4}_{1}(\tau) =
\begin{pmatrix}
i(\tau + \gr) & -\ti{w_{3}} \\
w_{3} & -i(\tau + \gr)
\end{pmatrix}, \;
A^{4}_{2} =
\begin{pmatrix}
1 & 0 \\
0 & -1
\end{pmatrix}.
\end{equation*}
We remark that system~\eqref{P4-JM} is related to the Jimbo-Miwa
system for $P_4$ (see equations (C.30) and (C.31) of
\cite{JM1981II}) by a simple gauge transformation and a change of
variables.

\subsection{Parametrization of Solutions in Terms of $P_4$}
 \label{subsec:P4-parametrization}

We begin with a parametrization of the general solution of
system~\eqref{P4sys} by the (general) solution of $P_4$. For this
purpose we use the correspondence between the $3\times3$ and
$2\times2$ Fuchs--Garnier representations for this system.

First of all we notice that system~\eqref{P4sys} admits first
integrals:
\begin{equation} \label{P4sysint}
w_{1}\ti{w_{1}} - w_{2}\ti{w_{2}} = 2i\gt_{0},\qquad w_{1}\ti{w_{1}}
+ w_{2}\ti{w_{2}} + w_{3}\ti{w_{3}} = 2i\gt_{\infty},
\end{equation}
where $\gt_{0}$ and $\gt_{\infty}$ are constants, which have a sense
of formal monodromies of the solution $\hat Y$ of
system~\eqref{P4-JM}. Motivated by the parametrization used by
Jimbo--Miwa (see equations \eqref{P4-JM} above) we define the
functions
\begin{equation} \label{P4-eqn}
\ti{y}(\tau)=-\frac{2w_1w_2}{\ti{w}_3}\qquad\mathrm{and}\qquad
\tilde z(\tau)=w_{1}\ti{w}_{1}.
\end{equation}
Now from \eqref{P4sys} we find that $\ti{y}$ and $\ti{z}$ satisfy
the following system of nonlinear ODEs:
\begin{equation}
 \label{eq:yz-P4}
\frac{d\ti{y}}{d\tau}=-4\ti{z}+\ti{y}^{2}+2i(\tau+\gr)\ti{y}+4i\gt_{0},\quad
\frac{d\ti{z}}{d\tau}=-\ti{y}\big(\ti{z}-i(\gt_{0}+\gt_{\infty})\big)-
\frac{2}{\ti{y}}\ti{z}(\ti{z}-2i\gt_{0}).
\end{equation}
Eliminating from this system the function $\tilde z(\tau)$ one finds
the following second order ODE for $\tilde y(\tau)$,
\begin{equation*}
\frac{d^{2}\ti{y}}{d\tau^{2}} = \frac{1}{2\ti{y}} \Big(
\frac{d\ti{y}}{d\tau} \Big)^{2} + \frac{3}{2} \ti{y}^{3} + 4i(\tau +
\gr) \ti{y}^{2} + 2\big(-(\tau + \gr)^{2} + i(1 - 2\gt_{\infty})
\big) \ti{y} + \frac{8\gt_{0}^{2}}{\ti{y}}.
\end{equation*}
We note that, under the change of variables $(\tau + \gr) \mapsto
e^{-i\pi/4}t$, $\ti{y} \mapsto e^{i\pi/4}y$, this equation is mapped
to the $P_{4}$ equation \eqref{P4} with $\ga = 2\gt_{\infty} - 1$
and $\gb = -8\gt_{0}^{2}$.

By using the parametrization for $\ti{y}(\tau)$ and $\ti{z}(\tau)$
given above we obtain the following expressions for the functions
$\{w_{j}(\tau),\ti{w}_{j}(\tau)\}$:
\begin{equation} \label{P4-w}
\begin{aligned}
w_{1} &= -\frac{f\ti{y}\ti{z}}{2},& \quad w_{2} &=
\frac{1}{g\ti{z}},& \quad
w_{3} &= -\frac{2g}{f} (\ti{z} - i\gt_{0} - i\gt_{\infty}), \\
\ti{w}_{1} &= -\frac{2}{f\ti{y}},& \quad \ti{w}_{2} &=
g\ti{z}(\ti{z} - 2i\gt_{0}),& \quad \ti{w}_{3} &= \frac{f}{g},
\end{aligned}
\end{equation}
where $f(\tau)$ and $g(\tau)$ satisfy the equations
\begin{equation} \label{P4-f,g}
\begin{aligned}
\frac{d}{d\tau}\log{f}&=-\frac{\ti{y}'}{\ti{y}}- \frac{1}{2}
\left( \frac{\ti{z}'}{\ti{z}} - \frac{\ti{y}(\ti{z} - i\gt_{0} - i\gt_{\infty})}{\ti{z}} +
\frac{2(\ti{z} - 2i\gt_{0})}{\ti{y}} \right), \\
\frac{d}{d\tau}\log{g}&=-\frac{\ti{z}'}{\ti{z}}- \frac{1}{2} \left(
\frac{\ti{z}'}{\ti{z} - 2i\gt_{0}} - \frac{\ti{y}(\ti{z} - i\gt_{0}
- i\gt_{\infty})}{\ti{z} - 2i\gt_{0}} + \frac{2\ti{z}}{\ti{y}}
\right).
\end{aligned}
\end{equation}
Integrating these equations we get:
\begin{equation} \label{P4-f,g-integrals}
\begin{aligned}
\frac{2}{f\ti{y}}&=(2\ti{z})^{1/2}\exp{\left(\int_{\tau_0}^{\tau}\Big(\frac{\ti{z}-2i\gt_0}
{\ti{y}}-\frac{\ti{y}}{2\ti{z}}(\ti{z}-i\gt_0-i\gt_{\infty})\Big)d\tau+i\ti{f}_0\right)},\\
\frac{1}{g\ti{z}}&=(2(\ti{z} -
2i\gt_{0}))^{1/2}\exp{\left(\int_{\tau_{0}}^{\tau} \Big(
-\frac{\ti{z}}{\ti{y}} + \frac{\ti{y}(\ti{z} - i\gt_{0} -
i\gt_{\infty})}{2(\ti{z} - 2i\gt_{0})} \Big) d\tau + i\ti{g}_{0}
\right)},
\end{aligned}
\end{equation}
where $\tau_{0}, \ti{f}_{0}, \ti{g}_{0}$ are constants of
integration. In fact, in general we need only two parameters, say
$f_0$ and $g_0$, while $t_0$ can be fixed. So, formulae
\eqref{eq:yz-P4}-\eqref{P4-f,g-integrals} solve the problem of
parametrization of general complex solutions of system \eqref{P4sys}
in terms of general complex solutions of $P_4$.

The formulae for the functions $f(\tau)$ and $g(\tau)$ can be
simplified by introducing the function $\ti{\gs}(\tau)$ following
the work of Jimbo--Miwa \cite{JM1981II}. We use $\ti{\gs}(\tau)$ to
eliminate the dependence on $\ti{y}$. The latter, as we see below,
is also important for specification of the physical solutions of the
3WRI system.

Consider the following identity, which can be proved just by
differentiation with the help of \eqref{P4sys},
\begin{equation*}
w_{1}w_{2}w_{3} + \ti{w}_{1}\ti{w}_{2}\ti{w}_{3} + i(\tau +
\gr)(w_{1}\ti{w_{1}} + w_{2}\ti{w_{2}}) = i\int_{\tau_{0}}^{\tau}
(w_{1}\ti{w_{1}} + w_{2}\ti{w_{2}}) d\tau.
\end{equation*}
Substituting in the expressions for $\{ w_{j}, \ti{w}_{j} \}$ given
in \eqref{P4-w} we find
\begin{equation*}
\ti{y}(\ti{z} - i\gt_{0} - i\gt_{\infty}) - \frac{2\ti{z}(\ti{z} -
2i\gt_{0})}{\ti{y}} = -2i(\tau + \gr)(\ti{z} - i\gt_{0}) +
2i\int_{\tau_{0}}^{\tau} (\ti{z} - i\gt_{0}) d\tau.
\end{equation*}
We define the function $\ti{\gs}(\tau)$ (cf. \cite{JM1981II}) by
\begin{equation} \label{P4-sigma}
i\ti{\gs}=-\ti{y}(\ti{z}-i\gt_{0}-i\gt_{\infty})+\frac{2\ti{z}(\ti{z}-2i\gt_{0})}{\ti{y}}-2i(\tau+\gr)\ti{z}.
\end{equation}
It follows from the above identity that we have for the derivative
of $\tilde\sigma(\tau)$,
\begin{equation}
 \label{eq:sigma-prime-P4}
\ti{\gs}' = -2\ti{z}.
\end{equation}
Summing up and subtracting definition~\eqref{P4-sigma} with the
second equation \eqref{eq:yz-P4} and solving the result with respect
to $\tilde y$ and $1/{\tilde y}$, respectively, one finds:
\begin{equation}
 \label{eq:P4-y-sigma}
\tilde y=-\frac{{\tilde z}'+i\tilde\sigma+2i(\tau+\rho)\tilde
z}{2(\tilde z-i\theta_0-i\theta_\infty)},\qquad \frac1{\tilde
y}=-\frac{{\tilde z}'-i\tilde\sigma-2i(\tau+\rho)\tilde z}{4\tilde
z(\tilde z-2i\theta_0)}.
\end{equation}
The compatibility condition of system~\eqref{eq:P4-y-sigma} with the
help of \eqref{eq:sigma-prime-P4} is equivalent to the following
SD-equation for the function $\tilde\sigma(\tau)$,
\begin{equation}
 \label{eq:sigma-P4}
\left(\frac{d^2\ti{\gs}}{d\tau^2}\right)^2=-4\left((\tau+\gr)\frac{d\ti{\gs}}{d\tau}-\ti{\gs}\right)^{2}-
4\frac{d\ti{\gs}}{d\tau}\left(\frac{d\ti{\gs}}{d\tau}+4i\gt_0\right)\left(\frac{d\ti{\gs}}{d\tau}+
2i\gt_0+2i\gt_\infty\right).
\end{equation}
With the help of the function $\ti{\gs}(\tau)$
equations~\eqref{P4-f,g-integrals} can rewritten as follows:
\begin{equation}
 \label{eq:P4-f,g-integrals-sigma}
\begin{aligned}
\frac2{f\ti{y}}&=\sqrt{2\ti{z}}\exp{\left(i\frac{(\tau+\rho)^2}2+
i\int_{\tau_{0}}^\tau\frac{\ti{\gs}}{2\ti{z}}d\tau + if_{0}\right)}, \\
\frac1{g\ti{z}}&=\sqrt{2(\ti{z}-2i\gt_{0})}\exp{\left(-i\frac{(\tau+\rho)^2}2-
i\int_{\tau_{0}}^{\tau}\frac{\ti{\gs}+4i\theta_0(\tau+\rho)}{2(\ti{z}-2i\gt_{0})}d\tau-ig_{0}\right)}.
\end{aligned}
\end{equation}
It is important to notice that formulae
\eqref{eq:P4-f,g-integrals-sigma}, \eqref{eq:P4-y-sigma}, and
\eqref{eq:sigma-prime-P4}, allow one to parameterize the functions
$w_j$ and $\tilde w_j$ for $j=1,2,3$ in terms of one function
$\tilde\sigma$. We call it $\sigma$-parametrization.

We are now ready to solve the reduced 3WRI system \eqref{RedSysIV}
in terms of $P_{4}$. Comparing the linear systems \eqref{LP-IVgen}
and \eqref{LP-P4} we get the following correspondence
\begin{subequations}
\begin{align}
v_{1}(\tau) &= 2iw_{1}(\tau),& v_{2}(\tau) &= -iw_{2}(\tau),&
v_{3}(\tau) &= -ie^{-i(\tau + \gr)^{2}} w_{3}(\tau),& \\
v^{*}_{1}(\tau) &= -i\ti{w}_{1}(\tau),& v^{*}_{2}(\tau) &=
2i\ti{w}_{2}(\tau),& v^{*}_{3}(\tau) &= -ie^{i(\tau + \gr)^{2}}
\ti{w}_{3}(\tau),&
\end{align}
\end{subequations}
where the functions $\{ w_{j}, \ti{w}_{j} \}$ are given in terms of
$\ti{y}$ and $\ti{z}$ by equations \eqref{P4-w} and
\eqref{P4-f,g-integrals}. This provide us the general similarity
solution for the coupled 3WRI system in terms of $P_4$.

In the physical situation we have to find the solution for which
$v_j=v_j^*$, where now the $*$ means complex conjugation. In this
case we impose the following restrictions on the parameters:
$$
\rho,t,t_0,f_0,g_0\in\mathbb R\qquad\mathrm{and}\qquad
\theta_0,\theta_\infty\in i\mathbb R.
$$
In that case the equation for $\tilde\sigma(\tau)$, equation
\eqref{eq:sigma-P4}, has real coefficients and we can take its
general real solution. Then using the $\sigma$-parametrization of
the functions $w_j$ and $\tilde w_j$ one proves that the physical
reduction is indeed fulfilled, provided $\tilde z>0$ and $\tilde z>2i\theta_0$.

\subsection{Relation between the Noumi--Yamada and Jimbo--Miwa\\
Fuchs--Garnier Pairs for the Fourth Painlev\'{e} Equation}
 \label{subsec:P4alt}

The symmetric form of $P_4$ \cite{B1992,VS1993,A1993,NY1999} reads:
\begin{equation}
 \label{eq:symmetric-P4}
\begin{aligned}
\frac{df_0}{dz}&=f_0(f_1-f_2)+\alpha_0,\\
\frac{df_1}{dz}&=f_1(f_2-f_0)+\alpha_1,\\
\frac{df_2}{dz}&=f_2(f_0-f_1)+\alpha_2,
\end{aligned}
\end{equation}
where $\alpha_k\in\mathbb C$ satisfy the following normalization condition,
$$
\alpha_0+\alpha_1+\alpha_2=1.
$$
The functions
\begin{equation}
 \label{eq:P4-symmetric-y}
y_k(t)=\sqrt{-2}f_k(z),\quad
t=z/\sqrt{-2},\qquad k=0,1,2
\end{equation}
solve $P_4$ \eqref{P4} for
\begin{equation}
 \label{eq:P4-symmetric-coeff}
\alpha=\alpha_{k+1(\mathrm{mod}\,3)}-\alpha_{k+2(\mathrm{mod}\,3)},\qquad
\beta=-2\alpha_k^2.
\end{equation}
The following Fuchs--Garnier pair in $3\times3$ matrices was obtained for \eqref{eq:symmetric-P4} by
Noumi and Yamada \cite{NY2000},
\begin{align}
 \label{NY-P4-FG-z}
\frac{d\Psi}{d\mu}&=-\left(
\begin{pmatrix}
0 & 0 & 0 \\
1 & 0 & 0 \\
f_0 & 1 & 0
\end{pmatrix}+\frac1\mu
\begin{pmatrix}
v_1 & f_1 & 1 \\
0 & v_2 & f_2 \\
0 & 0 & v_3
\end{pmatrix}
\right)\Psi\equiv-(A+\frac1\mu B)\Psi,\\
 \label{eq:NY-P4-FG-x}
\frac{d\Psi}{dz}&=-\left(
\mu\begin{pmatrix}
0 & 0 & 0 \\
0 & 0 & 0 \\
1 & 0 & 0
\end{pmatrix}+
\begin{pmatrix}
\frac z3-f_2 & 1 & 0 \\
0 & \frac z3-f_0 & 1 \\
0 & 0 & \frac z3-f_1
\end{pmatrix}
\right)\Psi=-(\mu P+Q)\Psi,
\end{align}
where the numbers $v_k$ and $\alpha_k$ are related as follows:
$$
\alpha_0=1+v_3-v_1,\qquad
\alpha_1=v_1-v_2,\qquad
\alpha_2=v_2-v_3.
$$
Note that, for a given set of the numbers $\{\alpha_k\}_{k=0,1,2}$, any one of the numbers $v_k$
can be taken arbitrarily. Using this fact we assume, without loss of generality, that
\begin{equation}
 \label{eq:P4-v1=1}
v_1=1.
\end{equation}
This assumption is equivalent to the additional transformation, $\Psi\rightarrow\mu^{1+v_1}\Psi$.

Now making for the solution $\Psi$ the generalized Laplace transform
$$
\Psi(\mu,z)=\int_Ce^{\mu \gz}\widetilde \Psi(\gz,z)d\gz
$$
with the suitably chosen contour $C$ (such that the corresponding terms appearing due to integration by parts
cancel) we arrive at the following system of equations for the Laplace image $\widetilde \Psi$:
\begin{equation}
 \label{eq:P4-FG-tilde-Psi}
\begin{aligned}
\frac{d\widetilde \Psi}{d\gz}&=(\gz I+A)^{-1}(B-I)\widetilde \Psi,\\
\frac{d\widetilde \Psi}{dz}&=\left(P(\gz I+A)^{-1}(B-I)-Q\right)\widetilde \Psi,
\end{aligned}
\end{equation}
In view of the condition \eqref{eq:P4-v1=1} the solution of
system~\eqref{eq:P4-FG-tilde-Psi} is reduced to the following system
in $2\times2$ matrices for the column vector
$\widetilde Y\equiv\widetilde Y(x,z)=(\widetilde \Psi_2(\gz,z),\widetilde \Psi_3(\gz,z))^T$, where
$x=1/\gz$ and $\widetilde \Psi_2,\widetilde \Psi_3$ are the
components of $\widetilde \Psi=(\widetilde \Psi_1,\widetilde
\Psi_2,\widetilde \Psi_3)^T$:
\begin{equation}
 \label{eq:P4-FG-tilde-Y}
\begin{aligned}
\frac{d\widetilde Y}{dx}&=\left(
\begin{pmatrix}
0 & 0 \\
-f_1 & -1
\end{pmatrix}x+
\begin{pmatrix}
f_1 & 1 \\
f_0f_1+v_2-1 & f_0+f_2
\end{pmatrix}+\frac1x
\begin{pmatrix}
1-v_2 & -f_2 \\
0 & 1-v_3
\end{pmatrix}
\right)\widetilde Y,\\
\frac{d\widetilde Y}{dz}&=\left(
\begin{pmatrix}
0 & 0 \\
-f_1 & -1
\end{pmatrix}x+
\begin{pmatrix}
f_0-\frac z3 & -1 \\
0 & f_1-\frac z3
\end{pmatrix}
\right)\widetilde Y.
\end{aligned}
\end{equation}
System \eqref{eq:P4-FG-tilde-Y} is mapped to the standard $2\times2$ Jimbo--Miwa Fuchs--Garnier pair for
$P_4$ (equations (C.30), (C.31)) of \cite{JM1981II}) by an appropriate triangular gauge transformation and
rescaling the variables $x \mapsto \sqrt{-2} x$, $z = \sqrt{-2}t$.  The resulting system is parameterized in
terms of the function $y_{2}(t)$, defined by \eqref{eq:P4-symmetric-y} for $k = 2$, which solves $P_4$ with
coefficients \eqref {eq:P4-symmetric-coeff}.

We conclude this subsection by noting that all of the transformations presented here and in
Subsection~\ref{subsec:P4-FG-Laplace} are invertible.  Thus it is straightforward to construct a direct mapping
between the $3\times3$ Fuchs--Garnier pair for $P_4$ obtained in this paper \eqref{LP-P4} and the one by Noumi and
Yamada for the symmetric form of $P_4$ \eqref{NY-P4-FG-z}.

\section{Similarity Reduction to the Third Painlev\'{e} Equation}
 \label{sec:P3}

The following similarity reduction was derived independently in
\cite{K1990} and \cite{MW1989}:
\begin{equation} \label{SimRedIII}
u_{1} = \exp[\tfrac{1}{2} x_{3} + \tfrac{i}{2} \gr x_{3}] v_{1}, \quad
u_{2} = \exp[\tfrac{1}{2} x_{3} + \tfrac{i}{2} \gr x_{3}] v_{2}, \quad
u_{3} = (x_{1} - x_{2})^{-1+i\gr} v_{3},
\end{equation}
where $v_j=v_j(\tau)$ with
\begin{equation}
\tau = (x_{1}-x_{2}) e^{x_{3}},
\end{equation}
and $\gr$ is a real constant. In this case system \eqref{3WRIcharac} reduces to the
system of ODEs:
\begin{equation} \label{RedSysIII}
\tau^{1+i\gr} v_{1}' = iv_{2}^{*}v_{3}^{*},\qquad
\tau^{1+i\gr} v_{2}' = -iv_{3}^{*}v_{1}^{*}\qquad
\tau^{i\gr} v_{3}' = iv_{1}^{*}v_{2}^{*}.
\end{equation}
It was shown in \cite{MW1989} that solutions of this system can be
represented in terms of an SD-type equation that is equivalent to
the particular case of the fifth Painlev\'{e} equation \eqref{P5}
with $\gd = 0$. It is well known \cite{G1975} that equation
\eqref{P5} with $\gd$ taken to be zero is equivalent to the general
$P_{3}$ equation.
\begin{remark}{\rm
As usual we generalize this similarity reduction to the coupled case of the 3WRI system.
To do it one adds to \eqref{SimRedIII} and \eqref{RedSysIII} the formally conjugated equations
\begin{equation*}
\begin{gathered}
u_1^*=\exp[\tfrac12x_3-\tfrac{i}2\gr x_3]v_1^*,\quad
u_2^*=\exp[\tfrac12x_3+\tfrac{i}2\gr x_3]v_2^*,\quad
u_3^*=(x_1-x_2)^{-1-i\gr}v_3,\\
\tau^{1-i\gr}{v_1^*}'=-iv_2v_3,\qquad
\tau^{1-i\gr}{v_2^*}'=+iv_3v_1,\qquad
\tau^{-i\gr}{v_3^*}'=-iv_1v_2.
\end{gathered}
\end{equation*}
Note that in the coupled case $\rho\in\mathbb C$ and the functions
$v_j$ and $v_j^*$ are not assumed to be complex conjugates. In the
most part of this Section we deal with the coupled 3WRI system and
turn back to the physical case at the end of
Subsection~\ref{subsec:P3-parametrization}.}
\end{remark}

\subsection{A $3\times3$ Fuchs--Garnier Pair for the Reduced System}

As in the case of $P_{4}$, $P_{5}$ and $P_{6}$ we will construct a
Fuchs--Garnier pair for the reduced system \eqref{RedSysIII} and
then carry out the explicit integration in terms of $P_3$. We
introduce the spectral parameter $\gl$ as
\begin{equation} \label{lambdaIII}
\gl = e^{-x_{3}}k.
\end{equation}
Taking $\kappa_{3}=0$ in \eqref{3WRI-LP} and writing $\Psi(x_{j},k)
= R(x_{j})\Phi(\tau,\gl)$ where $R(x_{j})$ is given by
\begin{equation*}
R(x_{1},x_{2},x_{3}) = \mathrm{diag}\, \Big( \exp[\tfrac{i}{2}\gr
x_{3}], \exp[-\tfrac{i}{2}\gr x_{3}], \exp[-\tfrac{1}{2}x_{3}]
\Big),
\end{equation*}
we obtain the following Fuchs--Garnier pair
\begin{subequations} \label{LP-IIIgen}
\begin{align}
M\Phi_{\gl} &= \Big( Q^{(1)} + \frac{Q^{(0)}}{\gl} \Big) \Phi \\
\Phi_{\tau} &= \Big( \gl P^{(1)} + P^{(0)} \Big) \Phi,
\end{align}
\end{subequations}
where the matrices $M, Q^{(1)}, P^{(1)}, Q^{(0)}, P^{(0)}$ are given
by
\begin{equation*}
M = \mathrm{diag}\, \big( 1, 1, 0 \big),
\end{equation*}
and
\begin{equation*}
Q^{(1)} = i\mathrm{diag}\, \big( -\tau\kappa_{2}, \tau\kappa_{1},
\kappa_{1} + \kappa_{2} \big), \quad P^{(1)} = i\mathrm{diag}\,
\big( -\kappa_{2}, \kappa_{1}, 0 \big),
\end{equation*}
\begin{equation*}
Q^{(0)} = i
\begin{pmatrix}
\tfrac{1}{2}\gr & \tau^{-i\gr}v_{3}^{*} & -v_{2} \\
\tau^{i\gr}v_{3} & -\tfrac{1}{2}\gr & v_{1}^{*} \\
-v_{2}^{*} & v_{1} & 0
\end{pmatrix}, \quad
P^{(0)} = i
\begin{pmatrix}
0 & \tau^{-1-i\gr}v_{3}^{*} & 0 \\
\tau^{-1+i\gr}v_{3} & 0 & 0 \\
-\tfrac12v_2^* & -\tfrac12v_1 & 0
\end{pmatrix}.
\end{equation*}

\subsection{Fuchs--Garnier Pairs for the Third Painlev\'{e} Equation}
 \label{subsec:P3-FG-Laplace}

As in the previous sections, to avoid a confusion with the
complex conjugation we introduce the following simpler notation:
\begin{subequations} \label{LP-P3}
\begin{align}
\begin{pmatrix}
1 & 0 & 0 \\
0 & 1 & 0 \\
0 & 0 & 0
\end{pmatrix}
\Phi_{\gl} &= \Big(
\begin{pmatrix}
\tau/2 & 0 & 0 \\
0 & -\tau/2 & 0 \\
0 & 0 & -1
\end{pmatrix}
+ \frac{1}{\gl}
\begin{pmatrix}
-\gt_{\infty}/2 & -\ti{w_{3}} & -w_{2} \\
w_{3} & \gt_{\infty}/2 & -\ti{w_{1}} \\
\ti{w_{2}} & w_{1} & 0
\end{pmatrix}
\Big) \Phi \label{LP3-1} \\
\Phi_{\tau} &= \Big(
\begin{pmatrix}
\gl/2 & 0 & 0 \\
0 & -\gl/2 & 0 \\
0 & 0 & 0
\end{pmatrix}
+ \frac{1}{\tau}
\begin{pmatrix}
0 & -\ti{w_{3}} & 0 \\
w_{3} & 0 & 0 \\
\tfrac\tau2\tilde w_2&-\tfrac\tau2w_1& 0
\end{pmatrix}
\Big) \Phi, \label{LP3-2}
\end{align}
\end{subequations}
where $\{w_{j}, \ti{w_{j}} \}$ are functions of $\tau$ and
$\gt_{\infty}$ is an arbitrary constant. The compatibility condition
for \eqref{LP-P3} is
\begin{equation*}
\begin{aligned}
\tau w_{1}' &= \ti{w_{2}}\ti{w_{3}},& \tau\ti{w_{1}}' &= w_{2}w_{3},& \\
\tau w_{2}' &= -\ti{w_{1}}\ti{w_{3}},& \tau\ti{w_{2}}' &= -w_{1}w_{3},& \\
\tau w_{3}' &= - \gt_{\infty}w_{3} + \tau\ti{w_{1}}\ti{w_{2}},&
\tau\ti{w_{3}}' &= \gt_{\infty}\ti{w_{3}} + \tau w_{1}w_{2}.&
\end{aligned}
\end{equation*}
We note that the third row of \eqref{LP3-1} gives the relation
\begin{equation*}
\gl \Phi_{3} = \ti{w_{2}} \Phi_{1} + w_{1} \Phi_{2},
\end{equation*}
and so we can eliminate $\Phi_{3}$ from the above system. The
resulting $2\times2$ system has the form:
\begin{subequations} \label{P3-JM}
\begin{align}
\frac{d\phi}{d\gl} &= \left( \frac{\tau}{2}
\begin{pmatrix}
1 & 0 \\
0 & -1
\end{pmatrix}
+ \frac{1}{\gl}
\begin{pmatrix}
-\gt_{\infty}/2 & -\ti{w_{3}} \\
w_{3} & \gt_{\infty}/2
\end{pmatrix}
- \frac{1}{\gl^{2}}
\begin{pmatrix}
w_{2}\ti{w_{2}} & w_{1}w_{2} \\
\ti{w_{1}}\ti{w_{2}} & w_{1}\ti{w_{1}}
\end{pmatrix} \right) \phi, \\
\frac{d\phi}{d\tau} &= \left( \frac{\gl}{2}
\begin{pmatrix}
1 & 0 \\
0 & -1
\end{pmatrix}
+ \frac{1}{\tau}
\begin{pmatrix}
0 & -\ti{w_{3}} \\
w_{3} & 0
\end{pmatrix} \right) \phi, \quad \phi = \begin{pmatrix} \Phi_{1} \\ \Phi_{2} \end{pmatrix}.
\end{align}
\end{subequations}
Making the change of variables $\tau \mapsto t^{2}, \gl \mapsto
\gl/t, w_{j}(\tau) \mapsto W_{j}(t)$, we find that system
\eqref{P3-JM} is equivalent (up to a diagonal gauge transformation)
to the Jimbo--Miwa system for the complete $P_{3}$ \cite {JM1981II}
in case $w_1\ti w_1+w_2\ti w_2=c_1\neq0$, or to the degenerate case
of $P_3$ (see e.g.~\cite{KV}) otherwise. Below we present the
corresponding parametrization for both cases.

\subsection{Parametrization of solutions in terms of $P_{3}$}
 \label{subsec:P3-parametrization}

Written in terms of the new variables the compatibility condition
for \eqref{LP-P3} becomes
\begin{equation} \label{P3-compatibility}
\begin{aligned}
tW_{1}' &= 2\ti{W_{2}}\ti{W_{3}},& t\ti{W_{1}}' &= 2W_{2}W_{3},& \\
tW_{2}' &= -2\ti{W_{1}}\ti{W_{3}},& t\ti{W_{2}}' &= -2W_{1}W_{3},& \\
tW_{3}' &= - 2\gt_{\infty}W_{3} + 2t^{2}\ti{W_{1}}\ti{W_{2}},&
t\ti{W_{3}}' & =2\gt_{\infty}\ti{W_{3}} + 2t^{2}W_{1}W_{2},&
\end{aligned}
\end{equation}
where $\{W_{j},\ti{W}_{j}\}$ are functions of $t$. This system
admits the following first integrals
\begin{equation}
 \label{eq:P3-first-integrals}
\begin{aligned}
&W_{1}\ti{W_{1}} + W_{2}\ti{W_{2}} = c_{1}, \\
&W_{1}W_{2}W_{3} - \ti{W_{1}}\ti{W_{2}}\ti{W_{3}} +
\frac{\gt_{\infty}}{2} \big( W_{1}\ti{W_{1}} - W_{2}\ti{W_{2}} \big)
= \frac{\gt_{0}}{2},
\end{aligned}
\end{equation}
where $c_{1}$ and $\theta_0$ are constants. Elementary computation
now shows that the function $y$ given by
\begin{equation} \label{y-P3}
y(t) = \frac{\ti{W_{3}}}{tW_{1}W_{2}},
\end{equation}
satisfies the $P_{3}$ equation \eqref{P3} with $\ga = 4\gt_{0}$,
$\gb=4(1-\gt_{\infty})$, $\gga = 4c_1^{2}$, $\gd = -4$. To
parameterize the functions $\{ W_{j}(t), \ti{W}_{j}(t) \}$ in terms
of $P_3$, we introduce the functions:
\begin{equation*}
z = tW_{1}\ti{W}_{1},\qquad w = W_{1}W_{2}.
\end{equation*}
Note that our notation for $w(t)$ is slightly different from the one
taken by Jimbo and Miwa \cite {JM1981II}. Using the expression for
$y$ given in \eqref{y-P3} and the compatibility conditions
\eqref{P3-compatibility}, we obtain the following system for $\{ y,
z, w \}$
\begin{equation*}
\begin{aligned}
t\frac{dy}{dt} &= 2(2z - c_1t) y^{2} + (2\gt_{\infty} - 1) y + 2t, \\
t\frac{dz}{dt} &= 4z(c_1t - z) y - (2\gt_{\infty} - 1) z + (\gt_{0} + c_1\gt_{\infty})t, \\
t\frac{d}{dt} (\mathrm{ln}\, w) &= 2(c_1t-2z)y.
\end{aligned}
\end{equation*}
It follows from \eqref{P3-compatibility} and the above expressions
that the functions $\{ W_{j}, \ti{W}_{j} \}$ are given as
\begin{equation} \label{eq:WJ}
\begin{aligned}
&W_{1}(t) = \frac{zf}{c_1t-z},\qquad
W_{2}(t) = g \frac{c_1t-z}{z}, \\
&\ti{W}_{1}(t) = \frac{c_1t-z}{tf},\qquad
\ti{W}_{2}(t) = \frac{z}{tg},\\
&W_{3}(t)=\frac{1}{tfg}\left( yz(c_1t-z)-\gt_{\infty}z+\frac{\gt_{0}+c_1\gt_{\infty}}{2}t\right),\\
&\ti{W}_{3}(t) = tfgy,
\end{aligned}
\end{equation}
where the functions $f(t)$ and $g(t)$ satisfy the following
equations
\begin{equation} \label{P3-f,g}
\begin{aligned}
t\frac{d}{dt} \log{f} &= 2y(c_1t - z) - t\frac{d}{dt} \log{\big(\frac{z}{c_1t-z}}\big), \\
t\frac{d}{dt} \log{g} &= -2yz - t\frac{d}{dt}
\log{\big(\frac{c_1t-z}{z}}\big).
\end{aligned}
\end{equation}
Now we are ready to solve the reduced 3WRI system \eqref{RedSysIII}
in terms of $P_3$. Just comparing systems \eqref{LP-IIIgen} and
\eqref{LP-P3} we find:
\begin{equation} \label{3WRI-P3}
\begin{aligned}
v_1(\tau)&=-iW_1(t),& v_2(\tau)&=-iW_2(t),& v_3(\tau)&=-it^{-2i\gr}W_3(t),& t&=\sqrt{\tau},\\
v^*_1(\tau)&=i\ti{W}_1(t),& v^*_2(\tau)&=i\ti{W}_2(t),&
v^*_3(\tau)&=it^{2i\gr}\ti{W}_3(t), &\gr&=i\gt_\infty,
\end{aligned}
\end{equation}
where the functions $W_j(t)$ are given by equations \eqref{eq:WJ}.
Equations~\eqref{3WRI-P3} provide a solution of the coupled 3WRI
system.

To get the solutions for the physical case we have to
guarantee that $v_{j}$ and $v^*_j$ are indeed complex conjugates. In
order to do this we notice that equations~\eqref{P3-f,g} imply the
following formulae for the functions $f$ and $g$:
\begin{align*}
f(t) &=
\frac{c_1t-z}z\sqrt{\frac{z}t}\,|t|^{\theta_\infty}\exp\left(-\frac{\theta_0+c_1\theta_\infty}2
\int^t_{t_0}\frac{dt}z+ic_2\right),\\
g(t)
&=\frac{z}{\sqrt{t}\sqrt{c_1t-z}}\,|t|^{\theta_\infty}\exp\left(\frac{\theta_0-c_1\theta_\infty}2
\int^t_{t_0}\frac{dt}{c_1t-z}+ic_3\right),
\end{align*}
where $t_{0}, c_{2}, c_{3} \in \mathbb{R}$ and the function $z(t)$
solves the following ODE,
\begin{multline}
 \label{eq:z-P3}
\frac{d^{2}z}{dt^{2}} = \frac{c_{1}t - 2z}{2z(c_{1}t - z)}\left(
\frac{dz}{dt} \right)^{2}+
\frac{z}{t(c_{1}t - z)} \left(\frac{dz}{dt}\right)+\frac{8z(c_1t-z)}t+\frac{c_1+4\theta_0\theta_\infty}{2t}\\
+\frac{(\theta_0-c_1\theta_\infty)^2-c_1^2}{2(c_1t-z)}-\frac{(\theta_0+c_1\theta_\infty)^2}{2z}.
\end{multline}
Therefore, if we choose $c_1\in\mathbb R$ and
$\theta_0,\theta_\infty\in i\mathbb R$, then the function $z(t)$ can
be taken real for real $t$ and should satisfy the following inequalities,
$0<z(t)/t<c_1$.
It is now easy to observe that under
these conditions equations~\eqref{3WRI-P3} define a similarity
solution for the physical reduction of 3WRI system.

In the case $c_{1} \neq 0$ the function $y(\tilde t)$ defined by the
following change of variables:
\begin{equation*}
z(t)=c_1t\frac{y(\tau)}{y(\tau)-1},\qquad \tau = t^2,
\end{equation*}
solves the degenerate case of the $P_5$ equation~\eqref{P5} with
the coefficients:
\begin{equation*}
\alpha=\frac{(\theta_0-c_1\theta_\infty)^2}{8c_1^2},\quad
\beta=-\frac{(\theta_0+c_1\theta_\infty)^2}{8c_1^2},\quad
\gamma=2c_1,\quad\delta=0.
\end{equation*}
In the special case $c_{1}=0$ equation~\eqref{eq:z-P3} coincides,
$z(t)=y(t)$, with the degenerate case of the $P_3$
equation~\eqref{P3} with the coefficients:
\begin{equation*}
\alpha=-8,\qquad\beta=2\gt_{0}\gt_\infty,\qquad\gamma=0,\qquad\delta=-\gt_0^2.
\end{equation*}

\subsection{Alternate Fuchs--Garnier Pairs for the Third Painlev\'{e} Equation}
 \label{subsec:P3alt}

To conclude this section we state without proof an alternate
reduction of system \eqref{LP-P3} to a $2\times2$ system. Using the
generalized Laplace transform \eqref{P6-LT} in \eqref{LP3-1}, the
resulting matrix equation has the form
\begin{equation}
\begin{pmatrix}
x - t/2 & 0 & 0 \\
0 & x + t/2 & 0 \\
0 & 0 & 1
\end{pmatrix}
\frac{d\ti{Y}}{dx} = -
\begin{pmatrix}
-\gt_{\infty}/2 +1 & -\ti{w_{3}} & -w_{2} \\
w_{3} & \gt_{\infty}/2 +1 & -\ti{w_{1}} \\
\ti{w_{2}} & w_{1} & 0
\end{pmatrix}
\ti{Y}.
\end{equation}
The determinant of the r.-h.s. matrix is zero by our choice of first
integrals for the $\{w_{j},\ti{w}_{j}\}$ system: by a special choice
of the parameter $c_1$ in \eqref{eq:P3-first-integrals}\footnote{Recall that $w_j(\tau)=W_j(t)$, see the end of
Subsection~\ref{subsec:P3-FG-Laplace}.}. We can therefore
make a gauge transformation $\ti{Y} = G \h{Y}$ where $G$ is the
diagonalizing matrix, to obtain
\begin{equation} \label{P5-delta=0}
\frac{d\h{Y}}{dx} = \left( \h{A}_{2} + \frac{1}{x - t/2} \h{A}_{1} +
\frac{1}{x + t/2} \h{A}_{0} \right) \h{Y},
\end{equation}
where the $\h{A}_{j}$ are all of the form
\begin{equation*}
\h{A}_{j} =
\begin{pmatrix}
* & * & 0 \\
* & * & 0 \\
* & * & 0
\end{pmatrix}.
\end{equation*}
Equation \eqref{P5-delta=0} can be
reduced (after a change of variables) to a $2\times2$ system of the
form
\begin{equation}
\frac{dY}{dx} = \left( A_{2} + \frac{1}{x-t} A_{1} + \frac{1}{x}
A_{0} \right) Y,
\end{equation}
where ${\rm tr} A_2=\det A_2=0$. Isomonodromy deformations in $t$ of
this equation are parameterized by solutions of the degenerate fifth
Painlev\'{e} equation \eqref{P5} with $\gd = 0$, see \cite{K2006},
\cite{KH1999}. It is interesting that in \cite{K2006} was found a
quadratic transformation relating isomonodromy deformations of
equations \eqref{P3-JM} and \eqref{P5-delta=0} for special values of
the corresponding formal monodromies. The results of this section
allow one to find a different transformation between these equations
and corresponding isomonodromy deformations for {\it all values} of
formal monodromies. We are going to present the details of this
correspondence in a separate work.

\bigskip \noindent
\textbf{Acknowledgements} This work was supported by ARC grant \#DP0559019.
The research was carried out during AVK's visits to the School of
Mathematics and Statistics at the University of Sydney, Australia.
We are also grateful to the referee for valuable comments which helped us to
improve the paper.

\appendix
\section{On Spectral Interpretation of the B\"{a}cklund Transformations for $P_5$}
 \label{app:A}

We recall the B\"{a}cklund transformation for $P_5$ that was found in
\cite{G1975} (see also~\cite{GL1982}):
\begin{align}
 \label{eq:backlund-y}
\hat y&=1-\frac{2\sqrt{-2\delta}ty}
{ty'-\sqrt{2\alpha}y^2+(\sqrt{2\alpha}-\sqrt{-2\beta}+t\sqrt{-2\delta})y+\sqrt{-2\beta}},\\
 \label{eq:backlund-alpha}
\sqrt{2\hat\alpha}&=\frac12\left(\frac\gamma{\sqrt{-2\delta}}+1-\sqrt{-2\beta}-\sqrt{2\alpha}\right),\\
 \label{eq:backlund-beta}
\sqrt{-2\hat\beta}&=\frac12\left(\frac\gamma{\sqrt{-2\delta}}-1+\sqrt{-2\beta}+\sqrt{2\alpha}\right),\\
 \label{eq:backlund-gamma}
\frac{\hat\gamma}{\sqrt{-2\hat\delta}}&=\sqrt{-2\beta}-\sqrt{2\alpha},\qquad
\sqrt{-2\hat\delta}=\sqrt{-2\delta}\neq0,
\end{align}
where $y=y(t)$ and $\hat y=\hat y(t)$ solve $P_{5}$ for the
parameters $\alpha$, $\beta$, $\gamma$, and $\delta$ and
$\hat\alpha$, $\hat\beta$, $\hat\gamma$, and $\hat\delta$,
respectively. The important feature of this transformation is that
the branches of the square roots in
equations~(\ref{eq:backlund-y})--(\ref{eq:backlund-gamma}) can be
taken arbitrarily but their choices remain the same in all the formulae.

Our goal here is to discuss the spectral interpretation of this
transformation. For this purpose we use the Jimbo-Miwa~\cite{JM1981II}
isomonodromy representation of $P_{5}$, the Fuchs--Garnier pair in
the Jimbo--Miwa parametrization.
Consider the following linear matrix ODE:
\begin{equation}
 \label{aeq:mainl}
\frac{d\Psi}{d\lambda} = \big(\frac t2\sigma_3 +
\frac{A_0}\lambda+\frac{A_1}{\lambda-1}\big)\Psi.
\end{equation}
Here
$\sigma_3=\left(\begin{array}{cc}
1&0\\0&-1
\end{array}\right)$
and the matrices $A_p$ ($p=0,1$) are independent of $\lambda$.
Assume the following parametrization of the matrices $A_p$,
\begin{equation}
 \label{aeq:JM-P5parametrization}
A_{0} =
\begin{pmatrix}
z + \frac{\gt_0}{2} & -u(z +\gt_{0}) \\
z/u & -z - \frac{\gt_0}{2}
\end{pmatrix},
\quad A_1 =
\begin{pmatrix}
-z -\frac{\gt_0+\gt_\infty}{2} & uy\big( z +\frac{\gt_0-\gt_1+\gt_\infty}{2} \big) \\
-\frac 1 {uy}\big( z+\frac{\gt_0+\gt_1+\gt_\infty}{2} \big) & z +
\frac{\gt_0+\gt_\infty}{2}
\end{pmatrix}.
\end{equation}
Then, the isomonodromy deformations of equation~(\ref{aeq:mainl}) with
respect to $t$ are governed by the following system of nonlinear
ODEs, which we will call the Isomonodromy Deformation System (IDS)
\begin{eqnarray}
t\frac{dy}{dt}&=&t
y-2z(y-1)^2-(y-1)\big(\frac{\gt_0-\gt_1+\gt_\infty}2y-\frac{3\gt_0+\gt_1+\gt_\infty}2\big),
 \label{aeq:ids1}\\
t\frac{dz}{dt}&=&yz\big(z+\frac{\gt_0-\gt_1+\gt_\infty}2\big)-\frac1y(z+\gt_0)
\big(z+\frac{\gt_0+\gt_1+\gt_\infty}2\big),
 \label{aeq:ids2}\\
t\frac
d{dt}\operatorname{log}u&=&-2z-\gt_0+y\big(z+\frac{\gt_0-\gt_1+\gt_\infty}2\big)+
\frac1y\big(z+\frac{\gt_0+\gt_1+\gt_\infty}2\big).
 \label{aeq:ids3}
\end{eqnarray}
In this system $\gt_\nu$ $(\nu=0,\,1,\,\infty)$ are complex
constants considered as parameters. Following \cite{JM1981II} we call them
the formal monodromies. Another widely used terminology for $\theta_\nu$ is 
the exponential differences, since they coincide with the differences of the
eigenvalues of the matrices
$A_\nu$\footnote{$A_\infty\equiv{\rm diag}(\, A_0+A_1)=-\theta_\infty\sigma_3/2$}.

Excluding the function $z$ from
equations~(\ref{aeq:ids1})--(\ref{aeq:ids2}) one finds that the function $y$
satisfies the fifth Painlev\'{e} equation (\ref{P5}) for the set of
the coefficients (\ref{eq:P5-coefficients-thetas}). In this case there exists a fundamental
solution of equation~\eqref{aeq:mainl} which solves the following equation:
\begin{equation}
 \label{aeq:maint}
\frac{d\Psi}{dt}=\left(\frac\lambda2\sigma_3 +
\frac1t(A_0+A_1+\frac{\theta_\infty}2\sigma_3)\right)\Psi.
\end{equation}

We note that by rescaling $t$ we may set the coefficients
$\hat\delta=\delta=-1/2$, and hence we may further put
\begin{equation*}
\sqrt{-2\delta}=\sqrt{-2\hat\delta}=\varepsilon=\pm1
\end{equation*}
in equations~(\ref{eq:backlund-y})--(\ref{eq:backlund-gamma}). To
take into account the possibility of different choices of branches
of the square roots in
equations~(\ref{eq:backlund-y})--(\ref{eq:backlund-gamma}) we
introduce the parameters $\varepsilon_1$, $\varepsilon_2$,
$\hat\varepsilon_1$, $\hat\varepsilon_1$, each taking the value
$\pm1$, in the following way
\begin{align}
 \label{eq:roots-theta}
\sqrt{2\alpha}&=\varepsilon_1\frac{\gt_0-\gt_1+\gt_\infty}2,\qquad
\sqrt{-2\beta}=\varepsilon_2\frac{\gt_0-\gt_1-\gt_\infty}2,\qquad
\gamma=1-\gt_0-\gt_1, \\
 \label{eq:roots-theta-^}
\sqrt{2\hat\alpha}&=\hat\varepsilon_1\frac{\hat\gt_0-\hat\gt_1+\hat\gt_\infty}2,\qquad
\sqrt{-2\hat\beta}=\hat\varepsilon_2\frac{\hat\gt_0-\hat\gt_1-\hat\gt_\infty}2,\qquad
\hat\gamma=1-\hat\gt_0-\hat\gt_1.
\end{align}
By substituting equations (\ref{eq:roots-theta}) and
(\ref{eq:roots-theta-^}) into
formulae~(\ref{eq:backlund-alpha})--(\ref{eq:backlund-gamma}) we get
the following equations relating the formal monodromies:
\begin{align}
 \label{eq:theta-epsilon-infty}
&\hat\gt_\infty=\varepsilon\frac{\hat\varepsilon_1-\hat\varepsilon_2}2(1-\gt_0-\gt_1)+
\frac{\hat\varepsilon_1+\hat\varepsilon_2}2\left(1-\frac{\varepsilon_1+\varepsilon_2}2(\gt_0-\gt_1)-
\frac{\varepsilon_1-\varepsilon_2}2\gt_\infty\right),\\
 \label{eq:theta-epsilon-0+1}
&\hat\gt_0-\hat\gt_1=\varepsilon\frac{\hat\varepsilon_1+\hat\varepsilon_2}2(1-\gt_0-\gt_1)+
\frac{\hat\varepsilon_1-\hat\varepsilon_2}2\left(1-\frac{\varepsilon_1+\varepsilon_2}2(\gt_0-\gt_1)-
\frac{\varepsilon_1-\varepsilon_2}2\gt_\infty\right),\\
 \label{eq:theta-epsilon-0-1}
&\hat\gt_0+\hat\gt_1=1+\varepsilon\left(\frac{\varepsilon_1-\varepsilon_2}2(\gt_0-\gt_1)+
\frac{\varepsilon_1+\varepsilon_2}2\gt_\infty \right).
\end{align}
Equations~(\ref{eq:theta-epsilon-infty})-(\ref{eq:theta-epsilon-0-1})
define $2^5=32$ different relations for the formal monodromies
$\hat\gt$'s, according to the number of tuples
$(\varepsilon,\varepsilon_1,\varepsilon_2,\hat\varepsilon_1,\hat\varepsilon_2)$.
It is easy to notice that all these formulae can be presented as the
compositions of the actions on the $\gt$-parameters of the
Schlesinger transformations ``dressing'' the points at infinity and zero
points:
\begin{align}
 \label{aeq:Spm}
S_{\pm,\pm}:\quad\gt_\infty\rightarrow\gt_\infty\pm1,\qquad
\gt_0\rightarrow\gt_0\pm1,\qquad\gt_1\rightarrow\gt_1,
\end{align}
with possibly the reflections:
\begin{align}
R_0:&\quad\gt_0\rightarrow-\gt_0,\qquad\gt_1\rightarrow\gt_1,
\qquad\gt_\infty\rightarrow\gt_\infty,\\
R_1:&\quad\gt_0\rightarrow\gt_0,\qquad\gt_1\rightarrow-\gt_1,
\qquad\gt_\infty\rightarrow\gt_\infty,\\
R_\infty:&\quad\gt_0\rightarrow\gt_0,\qquad\gt_1\rightarrow\gt_1,
\qquad\gt_\infty\rightarrow-\gt_\infty,\\
R_{01}:&\quad\gt_0\rightarrow\gt_1,\qquad\gt_1\rightarrow\gt_0,
\qquad\gt_\infty\rightarrow\gt_\infty,
\end{align}
and the following Okamoto-like transformation, mixing the
$\gt$-variables:
\begin{equation}
 \label{aeq:O}
{\cal O}:\quad\hat\gt_0=\frac{\gt_0+\gt_1-\gt_\infty}2,\quad
\hat\gt_1=\frac{\gt_0+\gt_1+\gt_\infty}2,\quad
\hat\gt_\infty=\gt_1-\gt_0.
\end{equation}
We call the last (mixing $\theta$'s) transformation the Okamoto transformation because
it coincides with the reflection $s_3={\cal O}$, introduced in Okamoto's work \cite{O1987b}.
In general, we call the Okamoto-like transformation any transformation for the Painlev\'e
equations that acts on the space of the corresponding formal monodromies as a linear operator with
the matrix having at least one row with nonzero elements: this condition means that it ``mixes''
all formal monodromies. Note that the absolute value of determinants of the matrices defining these
linear operators is always unity. These transformations appeared in Okamoto's studies of the Painlev\'e
equations ``on a regular footing" \cite{O1987a}--\cite{O1986}. They represent one of the reflections in the
subgroup of affine Weyl symmetries and exist for $P_4$, $P_5$, and $P_6$.

It is important to mention that the set of transformations~\eqref{aeq:Spm}--\eqref{aeq:O} is slightly wider
than those transformations that can be obtained via making compositions of B\"acklund
transformation~\eqref{eq:backlund-y} with different choices of the branches, e.g., from the latter transformation
one can obtain only the reflection $R_\infty\circ R_{01}$ rather than two of them separately. As we see below,
to produce  $R_\infty$ and $R_{01}$ we need an additional reflection of $t\to-t$, so that it is not a transformation
of the first kind in Okamoto's sense~\cite{O1987b}.

It is clear that among these transformations only 4 are independent, say,
$S_{+,-}$, $R_0$, $R_\infty$, ${\cal O}$, the others can be presented as follows:
$R_{01}=R_\infty\circ(R_\infty\circ{\cal O})^2$, $R_1=R_{01}\circ R_0\circ R_{01}$,
$S_{+,+}=R_{01}\circ S_{+,-}\circ R_{01}$, etc.

Now we are ready to discuss the presentation of these transformations on the solutions of the
Fuchs-Garnier pair \eqref{aeq:mainl}, \eqref{aeq:maint}.

It is well known that the Schlesinger
transformations $S_{\pm,\pm}$, are the gauge transformation of the $\Psi$ function,
$\Psi\to S(\lambda,t)/\sqrt{\lambda}\Psi$, where $S(\lambda,t)$ is a linear function of $\lambda$.
The general theory for these transformations in the framework of the isomonodromy deformations
can be found in \cite{JM1981II} and particular formulae for $P_5$ in \cite{MF1992}.

Reflections $R_0$ and $R_1$ do not have any spectral representation. For each $k=0,1$ both numbers:
$+\theta_k/2$ and $-\theta_k/2$, are the eigenvalues of the corresponding residue matrix $A_k$ and
because $\Psi$ is not normalized in the neighborhood of the singular points $\lambda=0$ and $\lambda=1$,
it does not ``feel a difference" between $+\theta_k$ and $-\theta_k$ $k=0,1$. The last statement means
that for each choice of the sign of $\theta_k$'s in the neighborhood of each singular point, $\lambda=0$
and $\lambda=1$, there exists its own series expansion representing the solution. However, this fact does
not affect the corresponding monodromy matrices or the residues. For example, consider $R_1$. Denote,
with hats the Painlev\'e functions after this reflection (see equations \eqref{aeq:mainl} and
\eqref{aeq:JM-P5parametrization}):
\begin{align}
R_1:&\qquad\hat\theta_0=\theta_0,\quad\hat\theta_1=-\theta_1,\quad\hat\theta_\infty=\theta_\infty
\nonumber\\
\text{from $A_0$:}&\qquad\hat z=z,\quad \hat u=u,
\nonumber\\
 \label{aeq:R_1}
\text{from $A_1$:}&\qquad\hat y\left(z+\frac{\theta_0+\theta_1+\theta_\infty}2\right)=
y\left(z+\frac{\theta_0-\theta_1+\theta_\infty}2\right).
\end{align}
The last equation represents the B\"acklund transformation of $P_5$ corresponding to the reflection $R_1$.
Substituting $z$ from equation~\eqref{aeq:ids1} one finds quite complicated explicit formula for $\hat y$
in terms of $y$. In an analogous way we find the nonlinear action of $R_0$:
\begin{align}
R_0:&\qquad\hat\theta_0=-\theta_0,\quad\hat\theta_1=\theta_1,\quad\hat\theta_\infty=\theta_\infty
\nonumber\\
\text{from $A_0$:}&\qquad\hat z-\frac{\theta_0}2=z+\frac{\theta_0}2,\quad\frac{\hat z}{\hat u}=\frac{z}{u},
\nonumber\\
 \label{aeq:R_0}
\text{from $A_1$:}&\qquad\hat u\hat y=uy\quad\Rightarrow \hat y=y\frac{z}{z+\theta_0},\quad\hat z=z+\theta_0,\quad
\hat u=u\frac{z+\theta_0}z.
\end{align}
Now we have two points of view: the Jimbo-Miwa parametrization is very good since it allows one to obtain
very easily quite nontrivial B\"acklund transformations for $P_5$; on the other hand, there is a discrepancy
between linear and nonlinear actions of the transformations $R_0$ and $R_1$, no linear action ``produces''
a nontrivial nonlinear action. Our result in Subsection~\ref{subsec:Okamoto-BT} related with the action of the
Okamoto transformation ${\cal O}$ means it is possible to obtain another parametrization of Fuchs-Garnier
pair \eqref{aeq:mainl}, \eqref{aeq:maint}, where $R_0$ and $R_1$ do not produce any nonlinear action on $y$,
see Remark~\ref{rem:P5Okamoto-parametrization}.

The situation with the reflection $R_\infty$ is different because of the normalization of the $\Psi$ function
at $\lambda=\infty$: the fact that at this point $\Psi$ has an irregular singularity does not play an essential
role in this question. The linear action reads
$$
\hat\Psi(\lambda,\hat t)=\sigma_1\Psi(\lambda,t)\sigma_1,\quad
\hat t=-t,\qquad
\sigma_1=\begin{pmatrix}
0&1\\
1&0
\end{pmatrix}
$$
The corresponding B\"acklund transformation is
$$
\hat y(\hat t)=1/y(t),\qquad
\hat z(\hat t)=-z(t)-\theta_0,\qquad
\hat u(\hat t)=1/u(t).
$$

The linear representation for $R_{01}$ is as follows:
$$
\hat\Psi(\hat\lambda,\hat t)=e^{-t\sigma_3/2}\Psi(\lambda,t),\qquad
\hat\lambda=1-\lambda,\qquad
\hat t=-t.
$$
It generates the following nonlinear representation,
$$
\hat y(\hat t)=1/y(t),\qquad
\hat z(\hat t)=-z(t)-\frac{\theta_0+\theta_1+\theta_\infty}2,\qquad
\hat u(\hat t)=e^{-t}y(t)u(t).
$$

The transformation that we obtain in Subsection~\ref{subsec:Okamoto-BT}, equation
\eqref{eq:P5Okamoto-transformations-yzu}, coincides with the composition $R_1\circ{\cal O}$. As is explained
above, $R_1$, being a nontrivial transformation for $P_5$, is not observable from ``spectral point of view",
so the same procedure gives us exactly ${\cal O}$ just by choosing a different parametrization of the
Fuchs-Garnier pair. The linear representation of the Okamoto transformation is given in \eqref{eq:P5Okamoto-integral}.

Let us present the ``true'' Jimbo--Miwa parametrization of the Fuchs-Garnier pair \eqref{aeq:mainl}, \eqref{aeq:maint}.
It is obtained from the original Jimbo--Miwa parametrization \eqref{aeq:JM-P5parametrization} by substituting in
it instead of $y$ and $z$ their expressions in terms of $\hat y$ and $\hat z$, obtained from
the first two equations in \eqref{eq:P5Okamoto-transformations-yzu}, a redefinition of
$u=\hat u\big(\hat z-(\theta_0-\theta_1+\theta_\infty)/2\big)$, the shift
$\hat z+\theta_1/2\to \hat z$, and removing the hats:
\begin{equation}
 \label{aeq:true-P5parametrization}
A_0=
\begin{pmatrix}
z-\displaystyle{\frac{\theta_\infty}2}&
-u\left(\!\left(z-\displaystyle{\frac{\theta_\infty}2}\right)^2-
\displaystyle{\frac{\theta_0^2}4}\right)\!\!\vspace{4pt}\\
\displaystyle{\frac1u}&
-z+\displaystyle{\frac{\theta_\infty}2}
\end{pmatrix}\!,\quad
A_1=
\begin{pmatrix}
-z&
uy\left(z^2-\displaystyle{\frac{\theta_1^2}4}\right)\!\!\vspace{4pt}\\
-\displaystyle{\frac1{uy}}&z
\end{pmatrix}\!.
\end{equation}
The system of isomonodromy deformations in this parametrization reads:
\begin{eqnarray}
t\frac{dy}{dt}\!\!\!&\!=\!&\!\!\!t y-2z(y-1)^2-\theta_\infty(y-1),
 \label{aeq:y-true}\\
t\frac{dz}{dt}\!\!\!&\!=\!&\!\!\!y\left(z^2-\frac{\theta_1^2}4\right)-\frac1y
\left(\!\left(z-\frac{\theta_\infty}2\right)^2-\frac{\theta_0^2}4\right),
 \label{aeq:z-true}\\
t\frac
d{dt}\operatorname{log}u\!\!\!&\!=\!\!\!&\!2\left(z-\frac{\theta_\infty}2\right)\left(\frac1y-1\right).
 \label{aeq:u-true}
\end{eqnarray}

Eliminating $z$ from equation~\eqref{aeq:y-true} and substituting it into equation~\eqref{aeq:z-true} one finds
that $y$ solves $P_5$ \eqref{P5} for the parameters:
\begin{equation}
 \label{aeq:P5-JM-true-coefficients}
\alpha=\frac{\theta_1^2}2,\qquad
\beta=-\frac{\theta_0^2}2,\qquad
\gamma=1+\theta_\infty,\qquad
\delta=-\frac12.
\end{equation}
Compare equations \eqref{aeq:y-true}--\eqref{aeq:u-true} and \eqref{aeq:P5-JM-true-coefficients} with the original
ones by Jimbo--Miwa \eqref{aeq:ids1}--\eqref{aeq:ids3} and \eqref{eq:P5-coefficients-thetas}.

Note, that now formulae \eqref{aeq:P5-JM-true-coefficients} look similar to the
analogous formulae for coefficients of the other Painlev\'e equations in the Jimbo--Miwa parameterizations.
To make this parametrization absolutely perfect we can apply  transformation $R_{01}$, to get $P_5$ with
$\alpha=\theta_0^2/2$ and $\beta=-\theta_1^2/2$.

We also remark that in this parametrization, changing of the signs of $\theta_0$ and $\theta_1$ has no effect
on $y$ and different transformations of $\Psi$ correspond to different transformations of $y$.
This was not the case for the original Jimbo--Miwa parametrization.

\end{document}